\input{style/arxiv-general.cfg}
\documentclass[MSNbibl,number,citesort,seceqn,dvips]{arxbj}
\makeatletter
   \@ifpackageloaded{graphicx}{}{\usepackage{graphicx}}
\makeatother
\usepackage{multirow,upgreek}
\usepackage{graphicx}
\usepackage{url,breakurl}

\aid{0}
\volume{21}
\issue{2}
\pubyear{2015}
\firstpage{1200}
\lastpage{1230}
\doi{10.3150/14-BEJ602} 

\makeatletter

\newcommand{\rrVert}{\Vert}
\newcommand{\rrvert}{\vert}
\newcommand{\llVert}{\Vert}
\newcommand{\llvert}{\vert}
\newtheorem{theorem}{Theorem}[section]
\newtheorem{prop}{Proposition}[section]
\newtheorem{lem}{Lemma}[section]
\newtheorem{cor}{Corollary}[section]
\newremark{rem}{Remark}[section]
\newremark{example}{Example}[section]
\newproclaim{Assumption}{Assumption}
\def\R{\mathbb{R}}
\def\E{\mathbb{E}}
\def\P{\mathbb{P}}
\def\d{\mathrm{d}}
\def\tr{\operatorname{tr}}
\def\cov{\operatorname{cov}}
\def\var{\operatorname{var}}
\def\bK{\mathbf{K}}

\newcommand{\trace}{\operatorname{trace}}

\newcommand{\diag}{\operatorname{diag}}

\renewcommand{\epsilon}{\varepsilon}

\def\vafrac#1#2{(#1)/(#2)}

\makeatother

\begin{document}
\begin{frontmatter}

\title{On the sample covariance matrix estimator of reduced effective
rank population matrices, with applications to fPCA}
\runtitle{The sample covariance matrix and reduced effective rank
population matrices}

\begin{aug}
\author[A]{\inits{F.}\fnms{Florentina} \snm{Bunea}\corref{}\thanksref{A}\ead[label=e2]{fb238@cornell.edu}} \and
\author[B]{\inits{L.}\fnms{Luo} \snm{Xiao}\thanksref{B}\ead[label=e1]{lxiao@jhsph.edu}}
\address[A]{Department of Statistical Science,
Cornell University, Ithaca, NY 14853, USA.\\
\printead{e2}}
\address[B]{Department of Biostatistics,
Bloomberg School of Public Health,
Johns Hopkins University,
Baltimore, MD 21205, USA.
\printead{e1}}
\end{aug}

\received{\smonth{12} \syear{2012}}
\revised{\smonth{1} \syear{2014}}

%
\begin{abstract}
This work provides a unified analysis of the properties of the sample
covariance matrix $\Sigma_n$ over the class of $p \times p$ population
covariance matrices $\Sigma$ of reduced effective rank $r_e(\Sigma)$.
This class includes scaled factor models and covariance matrices with
decaying spectrum.
We consider $r_e(\Sigma)$ as a measure of matrix complexity, and
obtain sharp minimax rates on the operator and Frobenius norm of
$\Sigma_n - \Sigma$, as a function of $r_e(\Sigma)$ and $\| \Sigma\|
_2$, the operator norm of $\Sigma$. With guidelines offered by the
optimal rates, we define classes of matrices of reduced effective rank
over which
$\Sigma_n$ is an accurate estimator. Within the framework of these
classes, we perform a detailed finite sample theoretical analysis of
the merits and limitations of the empirical scree plot procedure
routinely used in PCA. We show that identifying jumps in the empirical
spectrum that consistently estimate jumps in the spectrum of $\Sigma$ is
not necessarily informative for other goals, for instance for the
selection of those sample eigenvalues and eigenvectors that are
consistent estimates of their population counterparts. The scree plot
method can still be used for
selecting consistent eigenvalues, for appropriate threshold levels. We
provide a threshold construction and also give a rule for checking the
consistency of the corresponding sample eigenvectors. We specialize
these results and analysis to population covariance matrices with
polynomially decaying spectra, and extend it to covariance operators
with polynomially decaying spectra. An application to fPCA illustrates
how our results can be used in functional data analysis.
\end{abstract}

%
\begin{keyword}
\kwd{covariance matrix}
\kwd{eigenvalue}
\kwd{eigenvector}
\kwd{fPCA}
\kwd{high dimensions}
\kwd{minimax rate}
\kwd{optimal rate of convergence}
\kwd{PCA}
\kwd{scree plot}
\kwd{sparsity}
\end{keyword}

\end{frontmatter}

\section{Introduction}\label{sec:intro}

High dimensional covariance matrix estimation has received a high
amount of attention over the last few years.
This was largely motivated by the fact that the sample covariance
matrix $\Sigma_n$, based on a sample of size $n$, is not necessarily a
consistent estimator of the covariance matrix $\Sigma$ of a random
vector $X\in\R^p$, if $p>n$. In this regime, the shortcomings of
$\Sigma_n$ have been well understood for over a decade,
whenever we estimate a \textit{spiked} covariance matrix; see, for
instance, the seminal works of Baik and Silverstein \cite{Baik:06}
and Johnstone \cite{Johnstone:01}. By definition,
spiked models have a fixed number of
large eigenvalues and the rest equal to one. Therefore, the effective
number of parameters in such models is of order $p^2$, and there is no
hope to estimate them accurately from a small sample. To address this
issue, classes of \textit{sparse} covariance matrices have been
introduced in recent years. Depending on
the type of sparsity (entry-wise, row-wise, off-diagonal decay),
appropriate estimators have been introduced and shown to adapt to the
unknown sparsity structures, see, for instance, \cite{Bickel:08a,Bickel:08b,Cai:10,Cai:11}, among many others. It is important to note
that although sparse matrices, by definition, have a reduced number of
parameters, they can still be spiked. Therefore, the usage of the
sample covariance matrix $\Sigma_n$ in this context would still be
questionable, in addition to not rendering the appropriate sparse
structure. It is also of importance to observe that all sparse
covariance matrix models carry with them implicit modeling assumptions.
For instance, they are appropriate whenever many of the components of
$X$ are weakly correlated. They are also powerful for modeling
temporally or spatially
ordered variables, in cases where it is reasonable to assume that
variables apart in time or space have very little association.

However, there are many instances where these assumptions are not
satisfied, for example when the observed variables are known to have
strong associations with each other. If the association is
approximately linear, $\Sigma$ will be close to being a degenerate,
rank $r < p$ matrix, with possibly much fewer parameters than $p^2$, if
$r$ is small. To treat general, positive definite covariance matrices,
which have \textit{effectively} reduced rank, we
make use of the notion of \textit{effective rank}, first suggested by
Vershynin \cite{Vershynin:11} and given by
%
\begin{equation}
\label{re}r_e(\Sigma) = \frac{\trace(\Sigma)}{\|
\Sigma\|_2}.
\end{equation}

Here $\| \Sigma\|_2$ denotes the operator norm, or the largest
singular value, of $\Sigma$.
Clearly, $r_e(\Sigma)$ is smaller than the rank for degenerate
matrices and, in general, it can be significantly
smaller than $p$ if a large number of eigenvalues of $\Sigma$ are
relatively small.

Perhaps surprisingly, the finite sample properties of the sample
covariance matrix as an estimator of population matrices of reduced
effective rank are largely unstudied. For classes of matrices $\Sigma$
for which $r_e(\Sigma)$ and $\|\Sigma\|_2$ are appropriately bounded,
but allowed to vary with $n$ and $p$, we study the following problems:
\begin{enumerate}[(3)]
\item[(1)] Rate optimal estimation of $\Sigma$ via $\Sigma_n$, with respect
to the Frobenius and operator norms, in finite samples.

\item[(2)] Finite sample estimation of the location of a jump in the spectrum
of $\Sigma$, via $\Sigma_n$.

\item[(3)] Finite sample determination of the number of eigenvalues and
eigenvectors of $\Sigma_n$ that are accurate estimates, respectively,
of the eigenvalues and eigenvectors of $\Sigma$.

\item[(4)] Extensions of (2) and (3) to covariance operators, for functional
data.
\end{enumerate}

We study problem (1) in Section~\ref{sec:sample}. For data generated from a class of
sub-Gaussian distributions defined in Section~\ref{sec2.1}, we establish upper
bounds on the Frobenius norm $\|\Sigma_n - \Sigma\|_F$ and operator
norm $\|\Sigma_n - \Sigma\|_2$ that hold, with high probability, and
are near minimax optimal. We summarize these results in Table~\ref{table1}, which reveals that even if $p > n$, as long as
$r_e(\Sigma)$ and $\|\Sigma\|_2$ are appropriately small, $\Sigma_n$
continues to be an accurate estimator of $\Sigma$.
%
\begin{table}
\tablewidth=\textwidth
\tabcolsep=0pt
\caption{Optimal rates for the Frobenius and operator norm of $\Sigma
_n-\Sigma$: orders of magnitude depending on the regime of $p$. Within
each regime, the sizes of $r_e(\Sigma)$ and $\|\Sigma\|_2$ govern the rate}\label{table1}
\begin{tabular*}{\textwidth}{@{\extracolsep{\fill}}lll@{}}
\hline
Norm/values of $p$ & $p = \mathrm{O}(n^{\gamma}), \gamma\geq0$ &$p = \mathrm{O}\{\exp
(n)\}$ \\
\hline
Frobenius: $\|\Sigma_n-\Sigma\|_F$ & $\|\Sigma\|_2 \cdot r_e(\Sigma
)\cdot\sqrt{\frac{\ln n}{n}}$& $\|\Sigma\|_2\cdot r_e(\Sigma)\cdot
\sqrt{\frac{\ln n}{n}}$ \\
[3pt]
Operator: $\|\Sigma_n-\Sigma\|_2$&
$\|\Sigma\|_2 \cdot r_e(\Sigma)\cdot{\frac{\ln pn}{n}}$, if
$r_e(\Sigma) \geq\frac{n}{\ln pn}$&$\|\Sigma\|
_2\cdot r_e(\Sigma) \cdot\sqrt{\frac{\ln n}{n}}$
\\
& $ \|\Sigma\|_2 \cdot\sqrt{r_e(\Sigma)}\cdot\sqrt{\frac{\ln
pn}{n}}$, if $r_e(\Sigma) \leq\frac{ n}{\ln p n}$ &\\\hline
\end{tabular*}
\end{table}
The derivation of these bounds is presented in Section~\ref{sec:bounds}, where we
also study $\E\|\Sigma_n - \Sigma\|_F$ and $\E\|\Sigma_n - \Sigma
\|_2$, which have similar bounds, but sharper by $\ln n$ factors.
Guided by these results, we introduce and discuss classes of covariance
matrices of reduced effective rank, also in Section~\ref{sec:bounds}.



For problems (2) and (3), and their extension to (4), we investigate in
detail estimation performed by the ubiquitous scree plot method,
described below. Let $ \{\lambda_k, 1\leq k\leq p \}$,
arranged in descending order, denote the eigenvalues of $\Sigma$.
Similarly, let $\{\widehat{\lambda}_k, 1\leq k\leq p\}$, arranged in
descending order, denote the eigenvalues of the sample covariance
matrix $\Sigma_n$, henceforth called the sample eigenvalues. For a
given number $\tau$, called the threshold level, the scree plot method
consists in calculating the number
$K =: \max\{ k \dvtx  \widehat{\lambda}_k \geq\tau\}$
and retaining the $K$ largest sample eigenvalues. Typically, one also
retains the corresponding sample eigenvectors
$\widehat{\boldsymbol{\psi}}_k$, $ k \leq K$, for further analysis.
In Sections~\ref{sec:detectable} and \ref{sec:eigen}, we study when
this practice can be justified and for which threshold levels. To the
best of our knowledge, no theoretical study of the thresholding method
applied to $\Sigma_n$, of this nature, exists in the literature.

We study problem (2) in Section~\ref{sec:detectable}, where we give a
data-dependent construction of $\tau$ for
detecting minimal jumps in the spectrum of $\Sigma$. We say that a
minimal spectral jump occurs when there exists an index $s$ such that
$\lambda_s$ is a constant multiple of the noise level, and there is a
gap of at least the size of the noise level between $\lambda_s$ and
$\lambda_{s+1}$. The appropriate noise level for this class of
problems is proportional to $\E\|\Sigma_n - \Sigma\|_2$. The precise
definition and result are given in Theorem~\ref{t:detect}.
We apply this result to consistent estimation of the number of factors
in factor models in Example~\ref{example3}, complementing existing
methods, for example the AIC-type criterion in~\cite{baing}.

For population matrices with special structures, a spectral jump at the
minimal noise level may not exist. This is, for example, the case of
population matrices whose spectra exhibit a polynomial decay, which we
study in Section~\ref{sec:detectable:fpca}. In this case, spectral
jumps can still be detected, but they have to be larger than the noise
level, with order of magnitude depending on the rate of decay. We treat
this in Theorem~\ref{t:detect_fda}, where we offer guidance on a
data-dependent choice of $\tau$ for consistent jump detection under
this scenario.



We study problem (3) in Section~\ref{sec:eigen}. Finite sample bounds
on the difference between sample eigenvalues and eigenvectors and their
population counterparts
have been much less studied when $p > n$, and no unifying analysis over
the class of covariance matrices of reduced effective rank exists.
The study of consistent estimation of the eigenvalues and eigenvectors
of $\Sigma$ via $\Sigma_n$, in the classical asymptotic framework
where $p$ is fixed and $n \rightarrow\infty$, dates back half a
century, with notable works including those of Anderson \cite{Anderson:63} and
Muirhead \cite{muirhead}. Asymptotic analyses that allow
$p$ to grow with $n$
have been chiefly conducted in spike models, when $p/n$ converges to a
constant, and mostly concern the behavior of the largest sample
eigenvalue and corresponding eigenvector, see, for instance,  \cite{Johnstone:01} and  \cite
{Nadler}. None of these analyses can be
directly used or extended for studying problem (3). The most closely
related results to ours are those of
Kneip and Sarda \cite{Kneip:11}, who studied the finite sample
convergence rates of
the sample eigenvalues and eigenvectors of $\Sigma_n/p$ in factor
models, where $r_e(\Sigma/p)$ is finite and independent of $n$ and
$p$. We show in Section~\ref{sec:eigen} that their results are
particular cases of ours on studying problem (3) over classes of
population matrices of reduced effective rank. We show in Theorem~\ref
{t:eigenvalue_detect} that, for a given
desired precision level $\alpha$, we can construct a data-dependent
threshold level, which is a function of an estimate of the
minimum noise level and $\alpha$, such that all sample eigenvalues
above this threshold are close to the theoretical values at this
precision level, with high probability. A known result by Kneip and Utikal \cite{Kneip:01} can be used to show that, in general, it
would be misleading
to conclude that the sample eigenvectors corresponding to the sample
eigenvalues thus selected are also close to their population
counterparts. Our Theorem~\ref{t:eigenvector_detect} shows how to
complement the scree plot method by another simple strategy, in order
to further determine which sample eigenvectors are accurate estimates.
Interestingly, when the spectrum of $\Sigma$ decays polynomially, the
scree plot method once again suffices for
accurate estimation of both eigenvalues and eigenvectors and we make
this precise in Theorem~\ref{t:eigen_detect}.

In Section~\ref{sec:fpca}, we treat problem (4), by showing how the
results of the previous sections can be employed in fPCA. The data
consists in a sample of $n$ independent trajectories $X_i(t)$, of a
background stochastic process $X(t)$ with covariance operator $\mathcal
{K}$. Each trajectory
is observed at the same $m$ discrete data points $t_1 < t_2 < \cdots<
t_m$, and is corrupted by noise. Problem (2) has not been studied in
this context, but aspects of problem (3) have been thoroughly studied,
however only in asymptotic contexts. For perfectly observed
trajectories, at all time points $t$ and without additive noise,
Hall and Hosseini-Nassb \cite{Hall:06} use a result by Dauxois \textit{et al.} \cite{Dauxois:82} to develop a bootstrap
based approach for selecting the sample eigenvalues and eigenfunctions
that estimate the population counterparts at the parametric rate. For
discretely observed trajectories, the theoretical properties of the
estimated eigenvalues and eigenvectors have been established by, for
instance, Yao \textit{et al.} \cite{Yao:05}, Hall \textit{et al.}
\cite{Hall:06b} and Benko \textit{et al.} \cite{Benko:09}.
However, all these results are relative to the first few fixed eigenvalues
and eigenfunctions of $\mathcal{K}$, are of asymptotic nature, and the
selection of the appropriate number of features, in finite samples, is
left open. We bridge this gap here.

We study the class of covariance operators $\mathcal{K}$ with spectra
having polynomial decay, of which the Brownian motion is a chief
example. For this class, we show how the sample covariance matrix, in
connection with the scree plot method, can be employed to detect jumps
in the spectrum of the covariance operator, and to determine the number
of sample eigenvalues and eigenvectors that are accurate estimates of
the population eigenvalues and eigenfunctions, the latter evaluated at
the discrete observation points. Instrumental in this analysis, and new
relative to what we already developed in Sections~\ref{sec:detectable}
and \ref{sec:eigen}, are the results of Section~\ref{sec:fpca:finite}.

We denote by $\pi_m$ the projection mapping $X(t)$ into an
$m$-dimensional space $\R^m$, defined by $\pi_m(X) = (X(t_1), \ldots
, X(t_m))$. We refer to the distributions on $\R^m$ induced by $\pi
_m$ as the finite-dimensional distributions of $X$.
Let $\bK= m^{-1}\{\mathcal{K}(t_{j_1}, t_{j_2})\}_{1\leq j_1, j_2\leq
m}$ be the scaled covariance matrix corresponding to the
$m$-dimensional distribution of $X$. In Section~\ref{sec:fpca:finite},
we establish finite sample approximations of the eigenvalues and
eigenfunctions of the operator $\mathcal{K}$ by the eigenvalues and
eigenvectors of $\bK$. This allows us to transfer the assumptions on
the operator $\mathcal{K}$ to the matrix $\bK$, which in turn allows
us to apply the theory developed in Sections~\ref{sec:sample}--\ref
{sec:eigen} to functional data. Jump detection is presented in Section~\ref{sec:fpca:detectable} and the selection of the accurate sample
eigenvalues and eigenvectors is treated in Section~\ref{sec:fpca:eigen}.

The proofs of all our theoretical results are given in the \hyperref[sec:proof]{Appendix} and
in the supplemental material. We shall use the following notation
throughout our paper: $\|\cdot\|_F$, the Frobenius norm; $\|\cdot\|
_2$, the spectral/operator norm; $\|\cdot\|_1$, the nuclear norm; $\|
\cdot\| $, the Euclidean norm of a vector; $\tr(\cdot)$, the trace
of a square matrix; $I_p$, an identity matrix of dimension $p$. We will
also use the notation $\lesssim$ for inequalities that hold up to
multiplicative constants independent of $n$ and $p$ (or~$m$).
Throughout this paper, we regard a sample eigenvector $\widehat
{\boldsymbol{\psi}}$ as an estimate of its population counterpart
$\boldsymbol{\psi}$
and assume the sign of $\widehat{\boldsymbol{\psi}}$ is selected so
that $\widehat{\boldsymbol{\psi}}^{\prime}\boldsymbol{\psi}\geq0$.


\section{Some inequalities for the sample covariance matrix} \label
{sec:sample}
\subsection{Sub-Gaussian distributions}\label{sec2.1}
All the results of this paper are proved for a certain class of
sub-Gaussian distributions. In particular they all
hold for Gaussian vectors or processes. 
We recall that a zero-mean random variable $X\in\mathbb{R}$ is
\textit{sub-Gaussian} if there exists a constant $\sigma>0$ such that
$\E\exp(tX) \leq\exp(t^2\sigma^2/2)$, for all $t\in\mathbb{R}$.
Then it can be shown that $\sup_{k\geq1} k^{-1/2} (\E|X|^k)^{1/k}
<\infty$ and the sub-Gaussian\vspace*{1pt} norm of $X$ is defined to be $\|X\|
_{\psi_2} =\sup_{k\geq1} k^{-1/2} (\E|X|^k)^{1/k} $. A zero-mean
random vector $X\in\mathbb{R}^p$ is \textit{sub-Gaussian} if for any
non-random $u\in\mathbb{R}^p$, $u^{\prime}X$ is sub-Gaussian. The
sub-Gaussian norm of $X$ is defined as $\|X\|_{\psi_2} = \sup_{u\in
\mathbb{R}^p\backslash{\{0\}}} \|u^{\prime}X\|_{\psi_2}/ \|u\|$. We
will impose an additional assumption on a sub-Gaussian random vector:

\begin{Assumption}\label{as1} For a zero-mean sub-Gaussian random vector $X\in
\mathbb{R}^p$, we assume that there exists a constant $c_0>0$ such
that $\E(u^{\prime}X)^2 \geq c_0 \|u^{\prime}X\|_{\psi_2}^2$ for
all $u\in\mathbb{R}^p$.
\end{Assumption}

Assumption~\ref{as1} effectively bounds the higher moments of $X$ as polynomial
functions of the second moments of $X$. Let $\Sigma$ be the covariance
matrix of $X$, then $u^{\prime}\Sigma u \geq c_0 \|u^{\prime}X\|
_{\psi_2}^2$,\vspace*{1pt} for all $u\in\mathbb{R}^p$, under Assumption~\ref{as1}. We
will provide a number of distributions of interest that meet this
assumption below. Note first that if $X \in\mathbb{R}^p$ is
sub-Gaussian and satisfies Assumption~\ref{as1} and $O\in\mathbb{R}^{p\times
p}$ is an orthonormal matrix, then $OX$ is also sub-Gaussian and
satisfies Assumption~\ref{as1} with the same $c_0$.

%
\begin{example}
Let $ X \in\mathbb{R}^p$ be a random vector from a zero-mean
multivariate normal distribution. Then $X$ satisfies Assumption~\ref{as1} with
$c_0 = \uppi/2$ (\cite{Vu:12}).
\end{example}


\begin{example}
\label{example1}
Let $X=(X_1,\dots, X_p)^{\prime}$ and the components
$X_j$ are independent and have a zero-mean sub-Gaussian distribution.
Suppose there is a common constant $\sigma>0$ such that $\E\exp
(tX_j/\sqrt{\Sigma_{jj}} )\leq\exp(t^2\sigma^2/2)$ for all
$j$, where $\Sigma_{jj}$ is the variance of $X_j$. Then $X$ is
sub-Gaussian and satisfies Assumption~\ref{as1}. Moreover, if $\widetilde{X} =
OX$ where $O\in\mathbb{R}^{p\times p}$ is an orthonormal matrix, then
$\widetilde{X}$ is sub-Gaussian and satisfies Assumption~\ref{as1}.
\end{example}

A proof of the statements in Example~\ref{example1} is provided in
Appendix~\ref{sec:sample:example}.

\subsection{Accuracy of the sample covariance matrix over classes of
population matrices of reduced effective rank}
\label{sec:bounds}

Let $X_1, \ldots, X_n$ be i.i.d. observations of a random vector $X
\in\mathbb{R}^{p}$. Without loss of generality, we assume that $\E
(X) = 0$. Let $\bar{X} = n^{-1}\sum_{i=1}^n X_i$ and $\Sigma_n =
n^{-1}\sum_{i=1}^n (X_i-\bar{X}) (X_i-\bar{X})^{\prime}$ be the
sample covariance matrix. We establish below sharp probability upper
bounds on $\Sigma_n - \Sigma$, in terms of both the Frobenius and the
operator norms, as well as sharp bounds on the expectation of either
norm. The bounds stated below hold up to multiplicative constants
defined precisely in Appendix~\ref{sec:sample:proof}. Specifically,
$c_1, c_2$ and $c_3$ are defined in Propositions \ref{p:average}, \ref
{p:frobenius} and \ref{p:operator}, respectively. The constants are
independent of $n$ and $p$ and depend only on $c_0$ in Assumption~\ref{as1}. As
announced in the \hyperref[sec:intro]{Introduction}, we show that
the effective rank defined in (\ref{re}) governs the size of these
bounds. As a consequence, we introduce classes of population matrices
over which $\Sigma_n$ can be employed accurately even if $p > n$. In
some cases, we offer high-level practical guidance on assessing
whether, for a given data set, it is reasonable to assume the
covariance matrix of a generating distribution belongs to these classes.

%
\begin{theorem}
\label{t:frobenius}
Suppose $X$ is a random vector that satisfies Assumption~\ref{as1}. With
probability at least $1-5n^{-1}$,
\[
\|\Sigma_n-\Sigma\|_F \leq2c_1 \cdot\|
\Sigma\|_2 \cdot r_e(\Sigma ) \cdot\sqrt{
\frac{\ln n}{n}}.
\]
Furthermore,
\[
\E \bigl(\|\Sigma_n-\Sigma\|_F^2 \bigr)
\leq2\cdot\|\Sigma\| _2^2\cdot\frac{ r_e(\Sigma)^2}{n^2}\cdot \bigl
\{16c_1^2c_2 +1+2\bigl(c_1^2+c_1
\bigr)\exp(1) \bigr\}.
\]
\end{theorem}

%
\begin{theorem}
\label{t:operator}
Suppose $X$ is a random vector that satisfies Assumption~\ref{as1}.
With probability at least $ 1 - 4n^{-1}$,
\[
\|\Sigma_n -\Sigma\|_2 \leq(1+c_1+c_3)
\cdot\|\Sigma\|_2 \cdot \max \biggl\{\sqrt{\frac{r_e(\Sigma)\cdot\ln pn }{n}},
\frac
{r_e(\Sigma)\cdot\ln pn }{n} \biggr\}.
\]
Furthermore, with $C =: 2 \{5c_3^2 + 1+2(c_1^2+c_1)\exp(1)
\}$,
\[
\E \bigl(\|\Sigma_n-\Sigma\|_2^2 \bigr) \leq
C\cdot\|\Sigma\| _2^2 \cdot\max \biggl\{
\frac{r_e(\Sigma)\cdot\ln p }{n}, \biggl(\frac{r_e(\Sigma)\cdot\ln p }{n} \biggr)^2 \biggr\}.
\]
\end{theorem}

%
\begin{rem}
\begin{enumerate}[(iii)]
\item[(i)] As it can be seen from the proofs in
Appendix~\ref{sec:sample:proof}, all our results continue to hold if
$\Sigma$ is singular.

\item[(ii)] Probability bounds on $\|\Sigma_n - \Sigma
\|_2$, similar to those in Table~\ref{table1}, have been first derived
for distributions with bounded support in  \cite{Vershynin:11}, Section~5.4.3.

\item[(iii)] A probability bound on $\| \Sigma_n -
\Sigma\|_2$, of the same order of magnitude as the one given by
Theorem~\ref{t:operator}, has been established independently by
Lounici \cite{Lounici:13}, as this work developed.
However, our proof is based on a
version of Berstein's inequality for unbounded matrices, whereas
Lounici \cite{Lounici:13}
employs a version of this inequality developed for bounded matrices,
and therefore uses a different argument to complete his proof. The rest
of the results presented in Theorems \ref{t:frobenius} and \ref
{t:operator}, including the bounds on expected values in both cases
are, to the best of our knowledge, new. The rates given by Theorems
\ref{t:frobenius} and \ref{t:operator} above are minimax optimal
over the class of matrices with effective rank bounded by $\min(\sqrt {n}, p)$, up to logarithmic terms. We refer to Theorem~2 of
 \cite{Lounici:13} for the lower bound derivations with
respect to the
operator norm. The lower bound with respect to the
squared Frobenius norm derived in Theorem~2 of \cite{Lounici:13} is
of the order of $\|\Sigma\|_2^2\cdot r_e(\Sigma) \cdot p /n$ and is
larger than
the rate we derived in Theorem~\ref{t:frobenius}. However, the proof
of Theorem~2 in  \cite{Lounici:13} can be
tightened, by keeping only
the first line of his inequality (5.27), to show that the minimax lower
bound is in fact $\|\Sigma\|_2^2\cdot r_e^2(\Sigma) /n$. Therefore,
our rate is near minimax optimal, over the class of matrices of
effective ranks bounded by $\min (\sqrt{n}, p )$.

\item[(iv)] It is noteworthy that the sample estimator
$\Sigma_n$ is already minimax rate optimal, in both Frobenius and
operator norm, over the class of  matrices of effective ranks bounded
by $\min(\sqrt{n}, p)$. This suggests that, over this class, very
little can be gained from further thresholding or shrinking operations.
For instance, the nuclear norm penalized estimator, that would have
appeared to be a more appropriate estimator over this class,
has the same and optimal bound in operator norm (\cite{Lounici:13}),
and very similar performance to $\Sigma_n$ in the simulations we have
conducted.
\end{enumerate}
\end{rem}

In most situations, a scale-independent accuracy measure for $\Sigma
_n$ is desired. One such measure is provided
by the ratio between $\| \Sigma_n - \Sigma\|_F$ or $\|\Sigma_n -
\Sigma\|_2$ and $\|\Sigma\|_2$. Then, recalling that
$\lesssim$ denotes inequalities that hold up to multiplicative
constants, Theorems \ref{t:frobenius} and \ref{t:operator} show that,
with high probability,
%
\begin{equation}
\label{F} \frac{\|\Sigma_n - \Sigma\|_F}{\| \Sigma\|_2} \lesssim r_e(\Sigma )\sqrt{
\frac{\ln n}{n}},
\end{equation}
and
%
\begin{equation}
\label{oper} \frac{\|\Sigma_n - \Sigma\|_2}{\| \Sigma\|_2} \lesssim\max \biggl\{\sqrt{\frac{r_e(\Sigma)\cdot\ln pn }{n}},
\frac{r_e(\Sigma
)\cdot\ln pn }{n} \biggr\}.
\end{equation}
The above relative measures are informative even if $ \| \Sigma\|_2$
increases with $p$ and they motivate the introduction of the following
classes of population matrices. Let $\epsilon\in(0, 1) $ be a
complexity index that may decrease to zero with $n$ and $p$. Let
$\gamma\geq0$ be a given number. Define
\[
\mathcal{P}_{1}(\epsilon) := \biggl\{\Sigma\dvtx  r_e(\Sigma)
\lesssim \epsilon\frac{n}{\ln pn}; p = \mathrm{O}\bigl(n^{\gamma}\bigr) \biggr\},
\]
and
\[
\mathcal{P}_2(\epsilon) := \biggl\{\Sigma\dvtx  r_e(\Sigma)
\lesssim \epsilon\sqrt{ \frac{n}{\ln n}} \biggr\}.
\]

The definition of these classes resembles sparsity definitions in
sparse covariance matrix models, where a certain sparsity measure is
controlled. The introduction of $\mathcal{P}_1(\epsilon)$ or
$\mathcal{P}_{2}(\epsilon)$ complements therefore the literature on
sparse models, by advocating the study of low complexity models, where
$r_e(\Sigma)$ is used as a complexity measure. Then, similar to
existing results which show that accurate estimation over classes of
population covariance matrices of a certain low complexity level is
possible even if $p > n$, Theorems \ref{t:frobenius} and \ref
{t:operator} show that estimation of covariance matrices with reduced
effective ranks can also be performed accurately even if $ p > n$, as
long as the complexity index $\epsilon$ is appropriately small. And
this can be achieved, in terms of rate optimality, by the ubiquitously
used sample covariance matrix. Specifically:
\begin{enumerate}[(ii)]
\item[(i)] \textit{For any $n$ and $p$}, if $\Sigma\in\mathcal
{P}_2(\epsilon)$, then Theorems \ref{t:frobenius} and \ref
{t:operator} yield:
\[
\frac{\|\Sigma_n - \Sigma\|_2}{\| \Sigma\|_2} \leq\frac{\|\Sigma
_n - \Sigma\|_F}{\| \Sigma\|_2} \lesssim\epsilon,
\]
since $\|M\|_2 \leq\|M\|_F$ for any matrix $M$.
Thus, if $ \epsilon= \mathrm{o}(1)$, the scaled operator and Frobenius norms
will be small. Note that this size of $\epsilon$ implies that
$r_e(\Sigma) = \mathrm{o}(\sqrt{n/\ln n})$.

\item[(ii)] If $p = \mathrm{O}(n^{\gamma})$, $\gamma\geq0$, then Theorem~\ref{t:operator} guarantees the accuracy of $\Sigma_n$ with respect
to the operator norm over a larger class of population matrices, with a
less restrictive size of $r_e(\Sigma)$. Specifically, if
$\Sigma\in\mathcal{P}_1(\epsilon)$, then
\[
\frac{\|\Sigma_n - \Sigma\|_2}{\| \Sigma\|_2} \lesssim\epsilon,
\]
which is small as long as $\epsilon= \mathrm{o}(1)$, implying that $r_e(\Sigma
) = \mathrm{o} (n/\ln pn  )$.
We note that the restriction on the growth of $p$ is induced by the
explicit dependency on $p$ in the logarithmic term of the bound (\ref
{oper}), which makes this bound non-informative if $p = \mathrm{O} \{\exp
{(n)} \}$, or if $p \rightarrow\infty$ independently of $n$.
If this is the case, we can use the results from (i) above, which are
valid for any $n$ and $p$, albeit over a smaller class of population matrices.
\end{enumerate}

In general, it is challenging to determine whether the population
covariance matrix is in $\mathcal{P}_1(\epsilon)$ or $\mathcal
{P}_2(\epsilon)$, for some $\epsilon$. Whereas a full solution to
this problem is beyond the scope of this paper,
we offer guidance for a particular case below. It is based on the
following result, also independently derived by Lounici \cite{Lounici:13}.

%
\begin{theorem}
\label{t:trace}
For any random vector $X$ satisfying Assumption~\ref{as1},
\[
\bigl\llvert \tr(\Sigma_n)-\tr(\Sigma)\bigr\rrvert
\leq4c_1\sqrt{\frac
{\ln n}{n}}\cdot\tr(\Sigma),
\]
with probability at least $1-5n^{-1}$.
\end{theorem}

%

\begin{rem}
By Theorems \ref{t:operator} and \ref{t:trace} we have,
for any $p$ and $n$ large enough and with high probability that
\[
\|\Sigma_n - \Sigma\|_2 \leq\|\Sigma_n -
\Sigma\|_F \lesssim\tr (\Sigma) \sqrt{\frac{\ln n}{n}} \leq2\tr(
\Sigma_n) \sqrt{\frac
{\ln n}{n}} ,
\]
or
\[
\E\|\Sigma_n - \Sigma\|_2 \leq\E\|\Sigma_n -
\Sigma\|_F \lesssim\frac{\tr(\Sigma)} {\sqrt{n}} \leq\frac{2\tr(\Sigma
_n)} {\sqrt{n}}.
\]
Theorem~\ref{t:trace} provides direct practical guidance on the
accuracy of the un-scaled Frobenius and operator norm, irrespective of
the size of $\| \Sigma\|_2$. It shows that, as a first simple check,
one should compare $\tr(\Sigma_n)$ to $\sqrt{n}$ in order to decide
whether $\Sigma_n$ suffices as an estimator of $\Sigma$. This is
particularly useful when we have reasons to believe that the
population covariance matrix has a large number of very small eigenvalues.
\end{rem}

\begin{rem}
By Theorems \ref{t:operator} and \ref{t:trace}, we can
derive an upper bound for $r_e(\Sigma_n)$ as an estimator of
$r_e(\Sigma)$; see Theorem~\ref{t:effective} in Appendix~\ref{sec:thm_re}.
\end{rem}

\section{Detectable spectral jumps for population covariance matrices
of reduced effective rank}\label{sec:detectable}

In this section, we discuss consistent estimation of an index $s$ of a
population eigenvalue that is sufficiently separated from the next one,
and therefore sufficiently large itself. We will refer to such an index
as a jump. In what follows, sufficiently large and sufficiently
separated will be defined relative to the bounds
on $\E\|\Sigma_n-\Sigma\|_2$ given by Theorems \ref{t:frobenius}
and \ref{t:operator} in Section~\ref{sec:bounds}. We will use
a slightly enlarged, by a $\sqrt{\ln n}$ multiplicative factor,
version of these bounds, which yields the appropriate noise levels for
index consistency analysis, as illustrated in Theorem~\ref{t:correct}
below. Specifically, if $p = \mathrm{O}(n^{\gamma})$, for some $\gamma\geq0$,
the noise level is
%
\begin{equation}
\label{delta1}\eta_1 := C \|\Sigma\|_2\cdot\sqrt
{r_e(\Sigma)}\cdot\sqrt{\ln pn/n},
\end{equation}
and, if $p = \mathrm{O} \{\exp(n) \}$, the noise level is
%
\begin{equation}
\label{delta2} \eta_2 := C\|\Sigma\|_2\cdot
r_e(\Sigma)\cdot\sqrt{\ln n/n}.
\end{equation}
To avoid notational clutter, we introduced above a constant $C> 0$ to
bound all other constants appearing in the bounds of Theorems \ref
{t:frobenius} and \ref{t:operator}. Note that $C$ does not depend on
$n$ or $p$. For a data-dependent threshold $\widetilde{\tau}$, define
%
\begin{equation}
\label{jumpest} \widehat{s}(\widetilde{\tau}) := \max \{k \dvtx  \widehat{\lambda
}_k \geq\widetilde{\tau} \},
\end{equation}
where we recall that $\widehat{\lambda}_k$, $ 1 \leq k \leq p$, in
decreasing order, are the sample eigenvalues.
The following general theorem shows the interplay between
the quantities needed to define an index $s$ of the spectrum of $\Sigma
$ that can be regarded as a jump and consistently estimated and the
conditions required of a data-dependent thresholding level $\widetilde
{\tau}$ that makes $\widehat{s}(\widetilde{\tau})$ a consistent
estimate of $s$. Recall that ${\lambda}_k$, $ 1 \leq k \leq p$, in
decreasing order, are the population eigenvalues.

%
\begin{theorem}\label{t:correct} Let $j \in\{1, 2\}$ be fixed. Suppose
$\Sigma\in\mathcal{P}_{j}(\epsilon)$, for some $\epsilon\in(0,
1)$ and that Assumption~\ref{as1} holds. If there exist an index $s$ and
positive quantities $\tau_1$ and $\tau_2$ such that
%
\begin{equation}
\label{index} \lambda_{s} \geq\tau_1 +
\eta_{j} \quad \mbox{and}\quad  \lambda_{s+1} \leq\tau_2 -
\eta_{j},
\end{equation}
and a data-dependent threshold $\widetilde{\tau}$ that satisfies
%
\begin{equation}
\label{tau} 
\P ( \tau_2 \leq\widetilde{
\tau} \leq\tau_1 ) \geq1 - \delta
\end{equation}
for some $\delta\in(0, 1)$,
then
\[
\P \bigl(\widehat{s}(\widetilde{\tau}) = s \bigr) \geq1 - 5n^{-1} -
\delta.
\]
\end{theorem}

%
\begin{rem}
\begin{enumerate}[(ii)]
\item[(i)] Note that if (\ref{index}) holds, with
either $j = 1$ or $j = 2$, then implicitly
\[
\tau_1 \geq\tau_2 > \eta_{j}\quad  \mbox{and}\quad
\lambda_s - \lambda _{s+1} > 2 \eta_{j} + (
\tau_1 - \tau_2).
\]
Thus, condition (\ref{index}) re-states the well understood fact that
in order to estimate with high probability the index $s$ of what we
declare a large enough eigenvalue, at least larger than the noise
level, there must also be a gap larger than the noise level between
this eigenvalue and the one following it.

\item[(ii)] If an index $s$ satisfying (\ref{index})
exists, we will call it a jump in the spectrum of $\Sigma$ relative to
the triplet
$(\tau_1, \tau_2, \eta)$.
\end{enumerate}
\end{rem}


Theorem~\ref{t:correct} makes it clear that, for each $j \in\{1, 2\}
$, the minimal allowable values for $\tau_1$ and $\tau_2$
are of the order of $\eta_{j}$, and are larger than $\eta_{j}$. The
following result specializes Theorem~\ref{t:correct} to this situation
and offers a concrete construction of data-dependent thresholds that
satisfy (\ref{tau}) with $\delta= \mathrm{O}(n^{-1})$.
We begin by defining two data-dependent levels:
%
\begin{equation}
\label{polyn} \widetilde{\eta}_{1} = C\|\Sigma_n
\|_2 \cdot\sqrt{\frac
{r_e(\Sigma_n)\cdot\ln pn }{n}},
\end{equation}
and
%
\begin{equation}
\label{expo} \widetilde{\eta}_{2} = C\|\Sigma_n
\|_2\cdot r_e(\Sigma_n) \cdot \sqrt{
\frac{\ln n}{n}},
\end{equation}
where the constant $C$ is the same as in the definitions of $\eta_1$
and $\eta_2$.
We will also use the following notation throughout this section: we let
$c_1, c_2$ and $c_3$ be the constants defined in Section~\ref{sec:sample}, and we let
%
\begin{equation}
\label{const} \epsilon_1 = 4c_1 \sqrt{\ln n /n},\qquad
C_1 = 0.9,\qquad  C_2 = 1.
\end{equation}

%
\begin{theorem}\label{t:detect}
 Let $j \in\{1, 2\}$ be fixed. Suppose $\Sigma\in
\mathcal{P}_{j}(\epsilon)$, for some $\epsilon\in(0, 1)$ and that
Assumption~\ref{as1} holds. Let $\eta_j$ be defined by either (\ref{delta1})
or (\ref{delta2}) above. Assume that there exists an index $s_j$ such that
\[
\lambda_{s_{j}} \geq\frac{2(1+\epsilon_1)}{C_j}\eta_j +
\eta_j, \qquad \lambda_{s_{j} + 1} < 2C_j(1-
\epsilon_1)\eta_j-\eta_{j}.
\]
Then, if $j = 1$ and $(1+c_1 + c_3)\sqrt{\epsilon} < 0.19$,
\[
\P \bigl\{\widehat{s}(2\widetilde{\eta}_1) = s_{1} \bigr
\}\geq1-11n^{-1}.
\]
If $j=2$,
\[
\P \bigl\{\widehat{s}(2\widetilde{\eta}_2) = s_{2} \bigr
\}\geq1-11n^{-1}.
\]
\end{theorem}

%
\begin{rem}
Theorem~\ref{t:detect} shows that it is possible to
detect, with high probability, fine jumps, at the minimal level of the
noise levels quantified by (\ref{delta1}) or (\ref{delta2}),
respectively, via data-dependent thresholds. However, the expressions
given for $\widetilde{\eta}_1$ and $\widetilde{\eta}_2$ above
depend on unknown constants, that in turn depend on the unknown
distribution of the data.
For practical use, we suggest cross validation.
\end{rem}

\begin{example}[(Estimating the number of factors in a factor model)]
\label{example3}
Let $\Sigma$ be a covariance matrix arising from a finite factor model
(see, for example, \cite{Chamberlain:83,Fan:13}), with the
decomposition
%
\begin{equation}
\label{fact} \frac{\Sigma}{p} = \sum_{r=1}^{R}
\lambda_r \boldsymbol{\xi}_r \boldsymbol{
\xi}_r^{\prime} + \frac{\sigma^2}{p} I_p,
\end{equation}
where $R$ is a fixed number, $\lambda_1 >\cdots>\lambda_R >0$ can be
upper bounded independently of $p$, and $\boldsymbol{\Xi} =
[\boldsymbol{\xi}_1,\dots, \boldsymbol{\xi}_R]$ satisfies
$\boldsymbol{\Xi}^{\prime}\boldsymbol{\Xi} = I_R$. Then $\Sigma
/p$ has finite effective rank, $\eta_{2} = C\cdot\tr(\Sigma/p)
\sqrt{\ln n/n} = \mathrm{O} (\sqrt{\ln n/n} )$. Assume further that $n =
\mathrm{o}(p^2)$. Then, when $n$ is sufficiently large, both $ \sigma^2/p <
2(1-\epsilon_1)\eta_2-\eta_2$ and $\lambda_R + \sigma^2/p \geq
2(1+\epsilon_1)\eta_2 + \eta_2$ hold. By Theorem~\ref{t:detect},
$\widehat{s}(2\widetilde{\eta}_2/p)$ estimates $R$, the number of
factors, accurately, with high probability.
\end{example}

\subsection{Population covariance matrices with polynomially decaying
spectrum: Jump detection} \label{sec:detectable:fpca}
In this section, the analysis presented in Theorem~\ref{t:correct} is
specialized to
particular modeling assumptions on $\Sigma$. With a view towards
Section~\ref{sec:fpca}, in which we discuss functional data,
we treat here in more detail the class of population covariance
matrices whose spectrum exhibits a polynomial decay. Specifically, we
consider matrices satisfying the conditions below. Let $\mathit{EG}(\Sigma) :=
\{\lambda_1,\dots,\lambda_p\}$.

\begin{Assumption}\label{as2} There exist absolute constants $C_{1\lambda}$,
$C_{2\lambda}$ and $\beta_2\geq\beta_1>1$ such that $C_{2\lambda
}k^{-\beta_2}\leq\lambda_k \leq C_{1\lambda} k^{-\beta_1}$, for
all $k$. Moreover, there exist absolute constants $C_{3\lambda}$ and
$\beta_3>\beta_2$ such that $\min_{\lambda\in \mathit{EG}(\Sigma), \lambda
\neq\lambda_k}|\lambda-\lambda_k| \geq C_{3\lambda} k^{-\beta
_3}$, for all $k$.
\end{Assumption}

We will show in Section~\ref{sec:fpca} that these conditions appear
naturally in the study of data generated from the Brownian motion, and
in that case we give specific values for $\beta_1, \beta_2$ and
$\beta_3$.
Note that the largest eigenvalue of any population matrix $\Sigma$
satisfying Assumption~\ref{as2} is a constant independent of $p$. Moreover,
since $\beta_1 > 1$, the effective rank $r_e(\Sigma)$ of such
matrices will also have a constant value. Therefore, Assumption~\ref{as2}
ensures that $\Sigma$ belongs to $\mathcal{P}_2(\epsilon)$ with
$\epsilon\lesssim\sqrt{\ln n /n}$, irrespective of the value of $p$.
If $p = \mathrm{O}(n^{\gamma})$, then $\Sigma\in\mathcal{P}_2(\epsilon)$,
with $\epsilon\lesssim\ln pn /n$. Note further that the order of the
noise levels $\eta_1$ and $\eta_2$ given by (\ref{delta1}) and (\ref
{delta2}), respectively, are, under Assumption~\ref{as2}, $\sqrt{\ln pn/n}$
and $\sqrt{\ln n/n}$,
and therefore only differ by a $\sqrt{\ln p}$ factor when $r_e(\Sigma
)$ is a constant.

In the analysis below we will consider only $\eta_2 = \mathrm{O}(\sqrt{\ln
n/n})$, to allow for the possibility of $p$ growing independently of
$n$. We will specialize Theorem~\ref{t:correct} by determining the minimal
values of $\tau_1$ and $\tau_2$ that define a detectable jump. We
note that they will differ from the values given by Theorem~\ref
{t:detect}, which is not applicable to the class of models satisfying
Assumption~\ref{as2}. To see why, first notice that Theorem~\ref{t:detect}
presupposes the existence of an index $s$ such that $\lambda_s$ (or
$\lambda_{s+1}$) and $\lambda_s - \lambda_{s+1}$ are constant
multiples of the noise level $\eta_2$. An index with these properties
does not exist under Assumption~\ref{as2}. It is immediate to see why: assuming
that such an $s$ exists would imply that $\frac{1}{s^{\beta_{1}}} <
\frac{1}{s^{\beta_{3}}} $, which cannot hold for $\beta_3 > \beta_1$.
However, if the jump in the theoretical spectrum occurs at a level that
is larger, in order, than the noise level, then it can be again
detected, with high probability, as illustrated by the following theorem.

%
\begin{theorem}
\label{t:detect_fda}
Suppose $\Sigma$ satisfies Assumptions \ref{as1} and \ref{as2}. Assume $n$ is
sufficiently large such that the following mild technical condition holds,
\[
(1+\epsilon_1)^{\beta_1/\beta_3} - (1-\epsilon_1)^{\beta_1/\beta
_3}
< C_{1\lambda}^{-1} \bigl(3C_{3\lambda}^{-1}
\bigr)^{-\beta
_1/\beta_3}\eta_2^{1-\beta_1/\beta_3}.
\]
If there exists an index $s$ such that
\[
\lambda_s \geq \bigl\{C_{4\lambda}(1+\epsilon_1)
\eta_2 \bigr\} ^{\beta_1/\beta_3} + \eta_2,\qquad
\lambda_{s+ 1} < \bigl\{C_{4\lambda
}(1-\epsilon_1)
\eta_2 \bigr\}^{\beta_1/\beta_3} - \eta_2
\]
with $C_{4\lambda} = 3C_{2\lambda}^{-1}C_{1\lambda}^{\beta_3/\beta_1}$,
then
\[
\P \bigl\{ \widehat{s} \bigl((C_{4\lambda}\widetilde{\eta
}_2)^{\beta_1/\beta_3} \bigr) = s \bigr\} \geq1 - 11n^{-1}.
\]
\end{theorem}

%
\begin{rem}
\begin{enumerate}[(iii)]
\item[(i)] The technical condition holds for sufficiently large
$n$ because $(1+\epsilon_1)^{\beta_1/\beta_3} -\linebreak[4]  (1-\epsilon
_1)^{\beta_1/\beta_3} = \mathrm{O}(\epsilon_1) = \mathrm{O}(\eta_2) = \mathrm{o} (\eta
_2^{1-\beta_1/\beta_3} )$.

\item[(ii)] The discussion prior to Theorem~\ref{t:detect_fda}
above illustrates that attempting to determine spectral jumps in the
population matrix that occur at the minimal noise level may be an ill
posed problem for certain classes of covariance matrices. Theorem~\ref
{t:detect_fda} offers a solution when Assumption~\ref{as2} is met.

\item[(iii)] Under Assumption~\ref{as2}, Theorem~\ref{t:detect_fda}
shows that by setting $\widetilde{\tau} = (C_{4\lambda}\widetilde
{\eta}_2)^{\beta_1/\beta_3} =\linebreak[4]  \mathrm{O}_P \{ (\ln n/n)^{\beta
_1/(2\beta_3)} \}$ in (\ref{jumpest}) we can estimate the jump
with high probability.
\end{enumerate}
\end{rem}

\section{Accuracy of the sample eigenvalues and eigenvectors selected
via scree plot methods}\label{sec:eigen}
In this section, we investigate whether the eigenvalues and the
corresponding eigenvectors, obtained via the simple thresholding
method, for appropriate data-dependent thresholds, are accurate
estimates of their population counterparts.
We begin by discussing eigenvalue estimation. By Weyl's theorem, we
always have $\llvert \widehat{\lambda}_k - \lambda_k\rrvert  \leq\|
\Sigma_n - \Sigma\|_2$, for all $k$. However, this inequality is not
particularly informative when $\lambda_k$ is small, and the relative
difference $\widehat{\lambda}_k/\lambda_k - 1$ may be more
appropriate and is used here.
By combining Weyl's theorem and the results in Section~\ref{sec:bounds}, we obtain the following corollary.

\begin{cor}\label{c:eigenvalue}
Suppose that Assumption~\ref{as1} holds. Let $\eta_{\min}$ be either $\eta
_1$ or $\eta_2$, defined in (\ref{delta1}) and (\ref{delta2}).
\begin{enumerate}[(ii)]
\item[(i)] Then
%
\begin{equation}
\label{eigen_eq1} \biggl\llvert \frac{\widehat{\lambda
}_k}{ \lambda_k} - 1\biggr\rrvert \leq
\frac{\eta_{\min}}{\lambda_k},
\end{equation}
holds simultaneously for all $k$, with probability larger than $1 -
5n^{-1}$.

\item[(ii)] For any $n$ and $p$, and for all $k$, we have
%
\begin{equation}
\label{kneip} \frac{ |\widehat{\lambda}_k - \lambda
_k|}{p} \lesssim\frac{\tr(\Sigma)}{p}\sqrt{
\frac{\ln n}{n} },
\end{equation}
with probability larger than $1 - 5n^{-1}$.
\end{enumerate}
\end{cor}

\begin{example}[(Estimation of eigenvalues in a factor model)]
For the factor model (\ref{fact}) defined above,
Kneip and Sarda \cite{Kneip:11}, in their Theorem~2, bound the
left-hand side in (\ref
{kneip}), for all $k \leq R$, by a term of order $1/p + (\log
p/n)^{1/2}$, when $p > \sqrt{n}$. Under this scenario both their bound
and ours have the same order of magnitude, $\mathrm{O}(\sqrt{\ln n/n})$.
Corollary~\ref{c:eigenvalue} above shows that, moreover, the $\mathrm{O}(\sqrt {\ln n/n})$ rate
of convergence is still valid when (i) $p < \sqrt {n}$; (ii) $p$ grows independently of $n$ or $p = \mathrm{O}\{\exp(n)\}$, since
(\ref{kneip}) does not contain a $\log p$ factor.
\end{example}

Inequality (\ref{eigen_eq1}) of Corollary~\ref{c:eigenvalue} shows
that, in order to have $\llvert \widehat{\lambda}_k /\lambda_k -
1\rrvert \leq\alpha$, where $\alpha$ is a small number in $(0,1)$,
for all $k \leq K$, the index $K$ has to satisfy $\lambda_K \geq\eta
_{\min}/\alpha$. Note further that the last population eigenvalue
that can be accurately estimated only needs to be larger than this
threshold, and there are no further requirements on the relative size
of the following eigenvalue or on the size of their gap. Therefore,
taking $K$ equal to one of the estimators of the detectable jumps
derived in the previous section is unnecessary and would be misleading,
as in this way we would identify only the consistent estimates of those
population eigenvalues up to where jumps occur. The following theorem
shows how to identify the data-dependent number of sample eigenvalues
close to their population counterparts, under very mild assumptions.

%
\begin{theorem}\label{t:eigenvalue_detect}
Let $j \in\{1, 2\}$ be fixed. Suppose $\Sigma\in\mathcal
{P}_{j}(\epsilon)$, for some $\epsilon\in(0, 1)$ and that Assumption~\ref{as1} holds. For $\epsilon_1$ and $C_j$ defined in (\ref{const}) above,
and for some given $\alpha\in(0, 1)$, let
%
\begin{equation}
\label{econsist} \widetilde{K}_j = \max \biggl\{ k\dvtx  \widehat{
\lambda}_k \geq\frac{\widetilde{\eta
}_{j}}{C_j(1-\epsilon_1)}\bigl(1 + \alpha^{-1}\bigr)
\biggr\}
\end{equation}
for the data dependent $\widetilde{\eta}_j$ given by (\ref{polyn})
or (\ref{expo}) above.
Then, $\llvert \widehat{\lambda}_k /\lambda_k - 1\rrvert \leq\alpha
$, for all $k \leq\widetilde{K}_j$, with probability larger than $1-11n^{-1}$.
\end{theorem}

The study of the accuracy of the sample eigenvectors is more complex
and Proposition~\ref{p:eigenvector} below shows that the accuracy of
sample eigenvectors depends on both $\eta_{\min}$ and the gaps
between successive eigenvalues. Recall that $\boldsymbol{\psi}_k$ is
the eigenvector of $\Sigma$ associated with $\lambda_k$ and denote by
$\widehat{\boldsymbol{\psi}}_k$ the counterpart from the sample
covariance matrix. The sign of $\widehat{\boldsymbol{\psi}}_k$ is
selected so that $\widehat{\boldsymbol{\psi}}_k^{\prime}\boldsymbol
{\psi}_k\geq0$.

%
\begin{prop}
\label{p:eigenvector}
Let Assumption~\ref{as1} hold. Let $\eta_{\min}$ be given by either (\ref
{delta1}) or (\ref{delta2}). Assume that $\lambda_1 > \lambda_2 >
\cdots>\lambda_p>0$. Then, with probability $1 - 5n^{-1}$,
%
\begin{equation}
\label{eigen_eq4} \llVert \widehat{\boldsymbol{\psi}}_k- \boldsymbol{
\psi}_k\rrVert \leq\frac{\eta_{\min}}{\min_{\lambda\in \mathit{EG}(\Sigma), \lambda
\neq\lambda_k}|\lambda-\lambda_k|} +\frac{6(\eta_{\min})^2}{\min_{\lambda\in \mathit{EG}(\Sigma), \lambda\neq\lambda_k}|\lambda-\lambda_k|^2}
\end{equation}
for each $k = 1, \dots, n\wedge p$.
\end{prop}

The above proposition follows by combining Lemma A.1 in \cite{Kneip:01} with the results of Section~\ref{sec:bounds}.
Furthermore, by taking $\eta_{\min} = \eta_2$ and using the same
reasoning as the one following Corollary~\ref{c:eigenvalue}, we can
derive sharper bounds on the left-hand side of (\ref{eigen_eq4}) than
those derived, for factor models, in Theorem~2 of  \cite{Kneip:11}.
These bounds will hold for all $n$ and $p$, and are valid for more
general matrices $\Sigma$.

Proposition~\ref{p:eigenvector} makes it clear that, without further
information on the degree of separation of the spectrum of $\Sigma$,
the scree plot method applied to the spectrum of $\Sigma_n$, for any
data-adaptive threshold, cannot guarantee that the retained sample
eigenvectors are close to their population counterparts. The theorem
below provides a simple way for evaluating accuracy of estimated
eigenvectors, based on the gaps between consecutive sample eigenvalues.

%
\begin{theorem}\label{t:eigenvector_detect}
Let $0 < \alpha< 1$ be given and define $\widehat{\lambda}_0 =
+\infty$, and $\widehat{\lambda}_{p+1} = 0$. Let $j \in\{1, 2\}$ be
fixed. Suppose $\Sigma\in\mathcal{P}_{j}(\epsilon)$, for some
$\epsilon\in(0, 1)$ and that Assumption~\ref{as1} holds. Let $\epsilon_1$
and $C_j$ as defined in (\ref{const}) above, and let $\widetilde{\eta
}_j$ be given by (\ref{polyn}) or (\ref{expo}). Then for all $k\geq
1$ such that
%
\begin{equation}
\label{sample_spectrum} \min (\widehat{\lambda}_{k-1} - \widehat{
\lambda}_k, \widehat {\lambda}_k - \widehat{
\lambda}_{k+1} ) \geq\frac{\widetilde
{\eta}_{j}}{C_j (1-\epsilon_1)}\bigl(2 + 3\alpha^{-1}
\bigr),
\end{equation}
we have
$\llVert \widehat{\boldsymbol{\psi}}_k- \boldsymbol{\psi}_k\rrVert  \leq\alpha$, with probability larger than $1-11n^{-1}$.
\end{theorem}

%
\begin{rem}The theorem shows that, in order to establish
accuracy of a certain sample eigenvector, one just needs to check
whether (\ref{sample_spectrum}) holds. The procedure is general, but
$\widetilde{\eta}_j$ still depends on unknown constants that in turn
depend on the distribution of the data. We suggest the usage of a
cross-validation type criterion for practical use.
Also, note that if both consistent eigenvalue \textit{and} eigenvector
estimation is of interest, then one can use the scree plot method
outlined in Theorem~\ref{t:eigenvalue_detect} to determine the maximum
number of consistent eigenvalues, then use the procedure described in
Theorem~\ref{t:eigenvector_detect} to evaluate which of the retained
eigenvectors are also consistent.
\end{rem}

\subsection{Population covariance matrices with polynomially decaying
spectrum: Accurate feature estimation}\label{sec:eigen:fpca}

If Assumption~\ref{as2} holds, we have knowledge of the degree of separation of
the population spectrum. In this case, Theorem~\ref
{t:eigenvector_detect} suggests that we just need to find the largest
$k$ such that (\ref{sample_spectrum}) holds, since (\ref
{sample_spectrum}) will hold for all smaller $k$. Furthermore, under
Assumption~\ref{as2},
the two inequalities in (\ref{econsist}) and (\ref{sample_spectrum})
will be equivalent. This means that we can use again the scree plot
method and develop a data-dependent threshold $\widetilde{\eta}_{\mathrm{ev}}$
that guarantees both eigenvalue and eigenvector consistency. We
formalize this below, again using $\eta_2$ as the benchmark noise
level. Results in terms of $\eta_1$ can be derived in a similar manner.

%
\begin{theorem}\label{t:eigen_detect}
Let $0 <\alpha< 1$ be given. Suppose that Assumption~\ref{as1} holds and let
$\epsilon_1$ be given by (\ref{const}) above. Under Assumption~\ref{as2},
define
%
\begin{equation}
\label{vecs} \widetilde{\eta}_{\mathrm{ev}} = C_{1\lambda} \biggl[
\frac{3\widetilde{\eta}_2}{(1-\epsilon
_1)C_{3\lambda}\alpha} \biggr]^{\beta_1/\beta_3} +\frac{\widetilde
{\eta}_2}{1-\epsilon_1}
\end{equation}
for $\widetilde{\eta}_2$ given in (\ref{expo}). Let
\[
\widetilde{K}_{\mathrm{ev}} = \max \{k \dvtx  \widehat{\lambda}_{k}
\geq \widetilde{\eta}_{\mathrm{ev}} \}.
\]
Then $\llVert \widehat{\boldsymbol{\psi}}_k - \boldsymbol{\psi
}_k\rrVert  \leq\alpha$ and $\llvert \widehat{\lambda}_k/\lambda_k
-1 \rrvert  \leq\alpha/3$, for all $k\leq\widetilde{K}_{\mathrm{ev}}$, with
probability larger\vspace*{1pt} than $1 - 11n^{-1}$.
\end{theorem}


\section{An application to fPCA}\label{sec:fpca}

In this section, we specialize our results to the analysis of sample
covariance matrices constructed from functional data. Let $X_i(s),
i=1,\dots, n$, denote an i.i.d. sample of trajectories from a Gaussian
process $\{X(t)\dvtx  0\leq t\leq1\}$, with covariance function $\mathcal
{K}(s,t) = \cov\{X(s), X(t)\}$. Assume that we observe discretized
versions of these trajectories, possibly corrupted by noise
%
\begin{equation}
\label{data}Y_i(t_j) = \mu(t_j) +
X_i(t_j) + E_{ij},
\end{equation}
where $\mu(\cdot)$ is the mean function and $E_{ij}$ are mean zero
measurement errors that are independent of $X_i(\cdot)$. Assume $\var
(E_{ij})= \sigma^2$ is finite. We assume that all trajectories are
observed at the same set of $m$ points $\{0< t_1<t_2<\cdots
<t_{m-1}<t_m<1\}$ in $[0,1]$.

We consider classes of covariance operators satisfying the following
assumptions.

\renewcommand{\theAssumption}{\Alph{Assumption}}
\setcounter{Assumption}{0}
\begin{Assumption}\label{asA} $\mathcal{K}(s, t)$ is continuous and a
positive semi-definite kernel.
\end{Assumption}

Under Assumption \ref{asA}, Mercer's theorem guarantees that $\mathcal
{K}(s,t)$ admits the representation $\sum_{k=1}^{\infty} \rho_k \phi
_k(s)\phi_k(t)$,
where $\{\rho_1\geq\rho_2\geq\cdots\geq0\}$ are non-decreasing
eigenvalues and $\{\phi_k(\cdot), k=1,\dots\}$ are eigenfunctions
that are orthonormal in $L_2[0,1]$. Moreover,
$ \sum_k \rho_k =: \rho_0 <\infty$.

\begin{Assumption}\label{asB} $\sup_{k}\sup_{s\in[0,1]} \llvert \phi
_k(s)\rrvert $ is bounded by a constant $C_{5\lambda}$.
\end{Assumption}

\begin{Assumption}\label{asC} $\frac{\partial\mathcal{K}(t,t)}{\partial
t}$ is continuous in $(0, 1)$, right continuous at 0 and left
continuous at 1. Moreover, $\int\llvert \frac{\partial
K(t,t)}{\partial t}\rrvert  \,\mathrm{d}t$ is finite.
\end{Assumption}

\begin{Assumption}\label{asD} $\sup_{s\in[0,1]} \llvert \phi
_k^{(1)}(s)\rrvert \leq C_{6\lambda} k^{\gamma_1}$ for all $k$ where
$\phi_k^{(1)}(s)$ is the derivative of $\phi_k$ and $C_1, \gamma_1$
are positive constants. Here $\phi_k^{(1)}(0)$ is the right derivative
of $\phi_k$ at 0 and $\phi_k^{(1)}(1)$ is the left derivative of
$\phi_k$ at 1.
\end{Assumption}

Note that the trigonometric basis satisfies Assumptions \ref{asB}--\ref{asD}.

\begin{Assumption}\label{asE} Assumption~\ref{as2} of Section~\ref{sec:detectable:fpca} holds for the
eigenvalues of $\mathcal{K}$, and moreover, $\beta_1 >\gamma_1$.
\end{Assumption}

\begin{rem} All these assumptions hold for the Brownian
motion and the Brownian bridge on $[0,1]$, with $\beta_1 =\beta_2 = 2$,
$\beta_3 = 3$, and $\gamma_1 = 1$.
\end{rem}

We denote by $\pi_m$ the projection mapping $X(t)$ into an
$m$-dimensional space $\R^m$, defined by $\pi_m(X) = (X(t_1), \ldots
, X(t_m))$. We refer to the distributions on $\R^m$ induced by $\pi
_m$ as the finite-dimensional distributions of $X$.
Let $\bK= m^{-1}\{\mathcal{K}(t_{j_1}, t_{j_2})\}_{1\leq j_1, j_2\leq
m}$ be the scaled covariance matrix for the $m$-dimensional
distribution of $X$.
Let $ Y_i = \{Y_i(t_1),\dots, Y_i(t_m)\}^{\prime}$, $\bar{Y}(t) =
n^{-1}\sum_{i=1}^n Y_i(t)$ and
$\bar{Y} = \{\bar{Y}(t_1),\dots, \bar{Y}(t_m)\}^{\prime}$. To
facilitate the application of the results derived in the previous
sections to functional data we denote
%
\begin{equation}
\label{sigmafunc} \Sigma= \bK+ m^{-1}\sigma^2I_m.
\end{equation}
An estimate of $\Sigma$ is the scaled sample covariance matrix,
corresponding to discretely observed trajectories:
\[
\Sigma_n = m^{-1}n^{-1} \sum
_i (Y_i-\bar{Y}) (Y_i-
\bar{Y})^{\prime
}.
\]
To keep our presentation focused, we have employed the sample mean
$\bar{Y}$ as an estimator of the mean function of the process. For the scenario
we study below, of densely sampled trajectories, $\bar{Y}$ suffices.
One may also use a smooth estimator, but then an appropriate equivalent
of Proposition~\ref{p:average} will be needed and it is deferred to
future work.

We shall discuss in detail the usage of the scree plot method based on
the sample covariance matrix $\Sigma_n$ for estimating: (i) the
location of jumps in the spectrum of the covariance operator $\mathcal
{K}$; (ii) the number of sample eigenvalues and eigenvectors that are
accurate estimates of their population counterparts. The diagram below
illustrates the connections needed for this analysis.
\[
\bigl\{ \mathcal{K}(s,t)\bigr\}_{{s,t} \in[0, 1]} \longleftrightarrow_{1}
\bK =: m^{-1}\bigl\{\mathcal{K}(t_{j_1}, t_{j_2})
\bigr\}_{1\leq j_1, j_2\leq m} \longleftrightarrow_{2} \Sigma=: \bK+
m^{-1}\sigma^2I_m \longleftrightarrow_{3}
\Sigma_n.
\]
First, recall that we assumed (Assumption \ref{asE}) that the spectrum of the
covariance \textit{operator} $\mathcal{K}$ has polynomial decay, and
that in Sections~\ref{sec:detectable:fpca} and \ref{sec:eigen:fpca}
we addressed in detail (i) and (ii) for covariance \textit{matrices}
whose spectra have polynomial decay such that the largest eigenvalue
and the effective rank are both finite. In
order to employ these results here, we need to identify a matrix that
can be formed from $\mathcal{K}$ by evaluating it at a discrete set of
points and whose spectrum has essentially the same properties as that
of $\mathcal{K}$. For us, this matrix is $\mathbf{K}$ defined above:
without scaling it by $m$ their respective spectra cannot be close, as
they are not of the same order. We show this in Proposition~\ref
{p:approx} below and, moreover, we show that the eigenvectors of $\bK$
are close to the vectors formed by
evaluating the eigenfunctions of $\mathcal{K}$ at the time points
$(t_1, \ldots, t_m)$. Assumptions \ref{asB}--\ref{asD} above are crucial
for establishing these connections. To account for error terms in model
(\ref{data}), we will consider a slight modification of $\mathbf{K}$,
namely $\Sigma$ defined above in (\ref{sigmafunc}), which has
features similar to $\mathbf{K}$. We therefore expect that the scree
plot method applied to $\Sigma_n$ will lead to consistent estimates of
(i) and (ii) above, and we show that this is indeed the case in
Sections~\ref{sec:fpca:detectable} and \ref{sec:fpca:eigen}.


\subsection{Finite approximations of eigenfunctions and eigenvalues}
\label{sec:fpca:finite}
Here we provide a deterministic analysis of the quality of $\mathbf
{K}$ as an approximation to $\mathcal{K}$.
With slight abuse of notation, we denote the eigenvalues of $\bK$ by
$ \{\lambda_1, \lambda_2,\dots \}$ and the associated
eigenvectors by $ \{\boldsymbol{\psi}_1, \boldsymbol{\psi
}_2,\dots \}$. The eigenvalues and eigenvectors of $\Sigma$ are
then $\{\lambda_k + m^{-1}\sigma^2, \boldsymbol{\psi}_k\}$.
We let $\boldsymbol{\phi}_k = m^{-1/2}(\phi_k(t_1),\dots, \phi
_k(t_m))^{\prime}$. Note that we intend to compare $\boldsymbol{\psi
}_k$, which is an eigenvector and therefore has Euclidean norm equal to
one, with $\boldsymbol{\phi}_k$, hence the need for scaling in its
definition. We also denote by $\mathit{EG}(\mathcal{K})$ the spectrum of
$\mathcal{K}$. The following assumption is also needed.

\begin{Assumption}\label{asF} For the fixed design points $\{t_j\dvtx  1\leq
j\leq m\}$, there exists a constant $M>0$ such that $M^{-1}m^{-1} \leq
\min_{0\leq j\leq m} |t_{j+1}-t_j| \leq\max_{0\leq j\leq m}
|t_{j+1}-t_j| \leq Mm^{-1}$. Here $t_0 = 0, t_{m+1} = 1$.
\end{Assumption}

%
\begin{prop}\label{p:approx}
If Assumptions \textup{\ref{asA}}--\textup{\ref{asF}} hold and if $m$ is sufficiently large so that $
m^{(1-\beta_1)/(\beta_1+\gamma_1)} \leq1/12C_{7\lambda}$, for
$C_{7\lambda}$ given in (\textup{\ref{eq0}}), then we have
\label{prop3}
%
\begin{equation}
\label{eigenvalue} \sup_{k \geq1} \llvert \lambda_k -
\rho_k\rrvert \leq C_{8\lambda
} m^{\vafrac{1-\beta_1}{\beta_1+\gamma_1} },
\end{equation}
where $C_{8\lambda}= C_{5\lambda}^2C_{1\lambda}/(\beta_1-1)
+C_{1\lambda} + 13C_{7\lambda}\lambda_0$
and also
%
\begin{equation}
\label{trace} \bigl\llvert \tr(\bK)- \rho_0\bigr\rrvert \leq
C_{9\lambda} m^{-1}
\end{equation}
for some fixed positive constant $C_{9\lambda}$, independent of $m$.
Moreover, we have
%
\begin{eqnarray}
\label{eigenvector} \llVert \boldsymbol{\psi}_k - \boldsymbol{
\phi}_k\rrVert & \leq& \frac{C_{8\lambda} m^{(1-\beta_1)/(\beta_1+\gamma_1)}}{\min_{
\rho\in \mathit{EG}(\mathcal{K}),\rho\neq\rho_k}|\rho-\rho_k| }
\nonumber\\[-8pt]\\[-8pt]
&&{}+ 6 \biggl\{\frac{C_{8\lambda} m^{(1-\beta_1)/(\beta_1+\gamma
_1)}}{\min_{ \rho\in \mathit{EG}(\mathcal{K}),\rho\neq\lambda_k}|\rho
-\rho_k| } \biggr\}^2 + 7C_{7\lambda}m^{(1-\beta_1)/(\beta_1+\gamma
_1)}
\nonumber
\end{eqnarray}
for all $k\leq m^{1/(\beta_1+\gamma_1)}$.
\end{prop}

To the best of our knowledge, the result in Proposition~\ref{p:approx}
is new. The proof is given in Appendix~\ref{sec:finite:proof}. Whereas
the global result (\ref{trace}) is an immediate consequence of approximating
integrals by finite sums, the derivation of (\ref{eigenvalue}) and
(\ref{eigenvector}) is much more involved, and depends crucially on
the behavior of the spectrum and eigenfunctions of the covariance
operator $\mathcal{K}$. The combination of (\ref{eigenvalue}) and
(\ref{trace}) immediately yields the result below.

\begin{cor}\label{c:reduced_func}
Under the assumptions of Proposition~\ref{prop3}, $r_e(\bK) = \mathrm{O}(1)$
and, moreover, $r_e(\Sigma) = \mathrm{O}(1)$.
\end{cor}

This result shows that the finite dimensional distributions of
processes with eigenvalues decaying as in Assumption \ref{asE} automatically
have scaled covariance matrices $\mathbf{K}$ of finite effective rank.


\subsection{Detectable jumps in the spectrum of a covariance operator}
\label{sec:fpca:detectable}

The results derived in Theorem~\ref{t:detect_fda} can be easily
extended to the consistent estimation
of an index of the spectrum of the covariance operator where a jump
occurs. The following theorem shows that we can detect spectral jumps
of an appropriate size via a data driven thresholding of the spectrum
of $\Sigma_n$. Since Proposition~\ref{p:approx}
guarantees that the spectra of $\mathcal{K}$ and $\mathbf{K}$ are
close, the construction of these thresholds follows from Theorem~\ref
{t:detect_fda}.
%

%
\begin{theorem}
\label{t:detect_fda_f}
Suppose that $X(t)$ is a Gaussian process with a covariance function
that satisfies Assumptions \textup{\ref{asA}}--\textup{\ref{asF}}. The assumption on $m$ is the same as
in Proposition~\ref{p:approx}.
Let $\eta_2$ be given by (\ref{expo}). Assume $n$ is sufficiently
large such that the following mild technical condition holds,
\[
(1+\epsilon_1)^{\beta_1/\beta_3} - (1-\epsilon_1)^{\beta_1/\beta
_3}
< C_{1\lambda}^{-1} \bigl(3C_{3\lambda}^{-1}
\bigr)^{-\beta
_1/\beta_3}\eta_2^{1-\beta_1/\beta_3}.
\]
If there exists an index $s$ such that
\begin{eqnarray*}
\rho_s &\geq& \bigl\{C_{4\lambda}(1+\epsilon_1)
\eta_2 \bigr\} ^{\beta_1/\beta_3} + C_{8\lambda}m^{\vafrac{1-\beta_1}{\beta
_1+\gamma_1}} +
m^{-1}\sigma^2 + \eta_2,
\\
\rho_{s+ 1} &<& \bigl\{C_{4\lambda}(1-\epsilon_1)
\eta_2 \bigr\} ^{\beta_1/\beta_3} -C_{8\lambda} m^{\vafrac{1-\beta_1}{\beta
_1+\gamma_1}} -
m^{-1}\sigma^2 - \eta_2
\end{eqnarray*}
with $C_{4\lambda} = 3C_{3\lambda}^{-1}C_{1\lambda}^{\beta_3/\beta_1}$,
then
\[
\P \bigl\{ \widehat{s} \bigl((C_{4\lambda}\widetilde{\eta
}_2)^{\beta_1/\beta_3} \bigr) = s \bigr\} \geq1 - 11n^{-1}
\]
for $\widetilde{\eta}_2$ given by (\ref{expo}) above.
\end{theorem}

%
\begin{rem}  We have stated Theorem~\ref{t:detect_fda_f}
in terms of $\eta_2$ given by (\ref{delta2}) of Section~\ref{sec:detectable} above. Since $r_e(\Sigma)$ and $\|\Sigma\|_2$ are
finite in the context of Section~\ref{sec:eigen}, then $\eta_2 =
\mathrm{O}(\sqrt{\ln n/n})$. From the results of Section~\ref{sec:bounds},
summarized in Table~\ref{table1}, we recall that this is the optimal bound
on $\| \Sigma_n - \Sigma\|_2$, in the regime $m = \mathrm{O}\{\exp(n)\}$, as
$\eta_2$ is independent of $m$, and can therefore be employed even if
$m \rightarrow\infty$. This facilitates the direct translation of our
results to the ideal case of perfectly sampled trajectories, when $m =
\infty$. For each fixed $m$, the noise level $\eta_1$ given by (\ref
{delta1}), of order $\mathrm{O}(\sqrt{\ln nm/n})$ can also be employed, and in
this case the data adaptive threshold will be a function of $\widetilde
{\eta}_1$.
\end{rem}

\begin{rem}
Recall that for the Brownian motion $\beta_1 = \beta_2 = 2$, $\beta
_3 = 3$ and $\gamma_1 = 1$. In this case, Theorem~\ref
{t:detect_fda_f} shows that by thresholding the sample covariance
matrix at a level of $\mathrm{O}_P \{(\ln n/n)^{1/3} + m^{-1/3} \}$
we can identify the location of the population eigenvalue larger than
the minimal level $\mathrm{O} \{(\ln n/n)^{1/3} + m^{-1/3} \}$, as
long as the following eigenvalue is also $\mathrm{O} \{(\ln n/n)^{1/3} +
m^{-1/3} \}$ apart. This is similar to the results of Section~\ref{sec:detectable:fpca}. The difference resides in the existence of
the extra additive term $ m^{-1/3} $, which quantifies the
approximation error.
\end{rem}

\subsection{On the accuracy of the sample eigenvalues and eigenvectors
selected via thresholding methods
for functional data}\label{sec:fpca:eigen}

We specialize the results of Section~\ref{sec:eigen} for data
generated as in (\ref{data}), and when Assumptions\linebreak[4]  \textup{\ref{asA}}--\textup{\ref{asF}} hold. For
this, we first establish finite sample upper bounds for the sample
eigenvalues and eigenvectors.

\begin{prop}
\label{p:eigen_func} Suppose that $X(t)$ is a Gaussian process with a
covariance function that satisfies Assumptions \textup{\ref{asA}}--\textup{\ref{asF}}. The assumption on
$m$ is the same as in Proposition~\ref{p:approx}. Let $C_{10\lambda}
= \max(m^{-1}\sigma^2 + c_2\rho_0, C_{8\lambda})$ where $c_2$ is as
in Theorem~\ref{t:frobenius} and $C_{8\lambda}$ is as in Proposition~\ref{p:approx}. Define
%
\begin{equation}
\label{minf} \eta_f =: C_{10\lambda} \bigl(\eta_2
+ m^{\vafrac{1-\beta_1}{\beta
_1+\gamma_1}} \bigr).
\end{equation}
Then with probability at least $1-5n^{-1}$, the following holds for
each $k$:
\[
\llvert \widehat{\lambda}_k-\rho_k\rrvert \leq
\eta_f.
\]
Furthermore, with probability at least $1-5n^{-1}$, for each $1\leq
k\leq m^{1/(\beta_1+\gamma_1)}$,
%
\begin{eqnarray}
\label{eigenvector2} \llVert \widehat{\boldsymbol{\psi}}_k - \boldsymbol{
\phi}_k\rrVert &\leq&\frac{\eta_f}{\min_{\rho\in \mathit{EG}(\mathcal{K}), \lambda\neq
\rho_k}|\rho-\rho_k| } +\frac{6\eta_f^2}{\min_{\rho\in \mathit{EG}(\mathcal{K}), \rho\neq\rho
_k}|\rho-\rho_k|^2 }\nonumber\\[-8pt]\\[-8pt]
&&{} +
7C_{8\lambda}m^{\vafrac{1-\beta_1}{\beta
_1+\gamma_1}}.\nonumber
\end{eqnarray}
\end{prop}

The proof of Proposition~\ref{p:eigen_func} follows directly from
Proposition~\ref{p:approx}, Corollary~\ref{c:eigenvalue}, and
Proposition~\ref{p:eigenvector}, hence the details are omitted.

\begin{rem}
Proposition~\ref{p:eigen_func} evaluates the accuracy of sample
eigenvalues and eigenvectors as a function of both the sample size and
the number of observations per subject. In particular, for the Brownian
motion, we recall that $\eta_2 = \mathrm{O} \{(\ln n / n)^{1/2} \}$
and thus
\[
\llvert \widehat{\lambda}_k-\rho_k\rrvert \lesssim(\ln n
/ n)^{1/2} + m^{-1/3}\qquad  \mbox{for each } k
\]
with high probability. Reasoning as in Theorem~\ref
{t:eigenvalue_detect} of Section~\ref{sec:eigen}, it also follows that
the ratio between all sample eigenvalues above $\eta_{f}$, or above an
estimate of it, and the corresponding theoretical values, will also be
close to one, with high probability.
\end{rem}

We recall that the accuracy of the sample eigenvectors also depends on
how well separated the eigenvalues of the operator $\mathcal{K}$ are
from each other. Under our assumptions on the covariance operator, we
have control on the degree of separation. We can therefore derive the
analogue of Theorem~\ref{t:eigen_detect} of Section~\ref{sec:eigen}
for functional data, and state it below.

%
\begin{theorem}\label{t:eigen_func}
Assume the settings in Proposition~\ref{p:eigen_func} hold. Then, with
$\eta_{f}$ given by (\ref{minf}) above we define
\begin{eqnarray*}
{\eta}_{\mathrm{op}} = C_{1\lambda} \biggl(\frac{3\eta_f}{C_{3\lambda}\alpha}
\biggr)^{\beta_1/\beta_3} +\eta_f.
\end{eqnarray*}
Let
\[
{K}_{\mathrm{op}} = \max \{k \dvtx  \widehat{\lambda}_{k} \geq\eta
_{\mathrm{op}} \}.
\]
Then $\llVert \widehat{\boldsymbol{\psi}}_k - \boldsymbol{\phi
}_k\rrVert  \leq\alpha$, for all $k\leq\min\{{K}_{\mathrm{op}}, m^{1/(\beta
_1 + \gamma_1)} \}$, and $\llvert \widehat{\lambda}_k/\rho_k -1
\rrvert  \leq\alpha/3$, for all $k\leq{K}_{\mathrm{op}}$, with probability
larger than $1 - 11n^{-1}$.
\end{theorem}

%
\begin{rem}The proof is immediate, and identical to the
one of Theorem~\ref{t:eigen_detect} above. In light of Theorem~\ref
{t:eigen_detect}, the result above continues to hold when $\eta_{f}$
is replaced by an estimate; in order to keep the presentation clear we
contented ourselves here with the usage of the theoretical level $\eta
_f$. For the Brownian motion $\beta_1 = 2$, $\gamma_1 = 1$ and $\beta
_3 = 3$, resulting in
\[
\eta_f = \mathrm{O} \bigl\{ (\ln n / n)^{1/2} + m^{-1/3}
\bigr\} \quad \mbox{and} \quad {\eta}_{\mathrm{op}} = \mathrm{O} \bigl\{ (\ln n/n)^{1/3} +
m^{-2/9} \bigr\}.
\]
Reasoning as in Section~\ref{sec:eigen}, we conclude that a
thresholding level that is larger than the minimal ${\eta}_{\mathrm{op}}$ guarantees
the accuracy of the sample eigenvalues and eigenvectors. For the
Brownian motion, the number of accurate sample eigenvectors is always
upper-bounded by $m^{1/3}$, but it may be smaller, depending on the
relative value of $K_{\mathrm{op}}$.
\end{rem}

\begin{appendix}
\section{Technical proofs}\label{sec:proof}
The proofs for the lemmas, propositions and theorems not included below
are provided in the supplemental article (\cite{Bunea:14b}).
\subsection{Technical proofs of Section \texorpdfstring{\protect\ref{sec:sample}}{2}}
\subsubsection{Three useful lemmas}

%
\begin{lem}
\label{l:sub1}
Let $X\in\mathbb{R}^p$ be a generic vector. Let $\Delta= \{u = (u_1,
\dots, u_p)^{\prime}\in\mathbb{R}^p\dvtx  |u_1| = \cdots= |u_p| = 1\}$.
Then for any positive integer $d$,
\[
\|X\|^{2d} \leq\frac{1}{2^p} \sum_{u\in\Delta}
\bigl(u^{\prime}X\bigr)^{2d}.
\]
\end{lem}

\begin{rem} In the following proofs, we will assume
sometimes, without loss of generality, that $\Sigma$ is a diagonal
matrix. This can be immediately justified as follows. Consider the
eigendecomposition $\Sigma= ODO'$, where $O$ is an orthonormal matrix
and $D$ is a diagonal matrix. Then $\cov(O'X) = D$ and $\|X\| = \|O'X\|
$. Similar arguments can be employed when we consider orthonormal
transforms of matrices, and evaluate either their Frobenius or operator norm.
\end{rem}

\begin{lem}
\label{l:sub2} Let $X\in\mathbb{R}^p$ be a zero-mean sub-Gaussian
random vector that satisfies Assumption~\ref{as1}. For any positive integer $d$,
\[
\E\|X\|^{2d} \leq\frac{(2d)^d}{c_0^d} \bigl[\tr(\Sigma)
\bigr]^d.
\]
\end{lem}

\begin{lem}
\label{l:sub_norm}
Let $X\in\mathbb{R}^p$ be a zero-mean sub-Gaussian random vector and
satisfies Assumption~\ref{as1}. Then
\[
\|\|X\|\|_{\psi_2}^2 \leq\frac{2\tr(\Sigma)}{c_0}.
\]
\end{lem}

\subsubsection{Proof of the statements in Example \texorpdfstring{\protect\ref{example1}}{2.2}}\label{sec:sample:example}\vspace*{-12pt}
\begin{pf*}{Proof of the statements in Example~\ref{example1}} We only
need to show that $X$ is sub-Gaussian and satisfies Assumption~\ref{as1}. Let
$u\in\mathbb{R}^p$ be an arbitrary non-random vector. Then for any
$t\geq0$,
\[
\E\exp\bigl(tu^{\prime}X\bigr) = \prod_{j=1}^p
\E\exp(tu_j X_j) \leq\prod
_{j=1}^p \exp \bigl\{ (tu_j\sqrt{
\Sigma_{jj}} )^2\sigma ^2/2 \bigr\}=\exp \bigl
\{t^2 \bigl(u^{\prime}\Sigma u\bigr)\sigma^2/2 \bigr
\},
\]
where the last equality holds because $\Sigma$ is a diagonal matrix as
the components of $X$ are independent. Hence, $u^{\prime} X$ is
sub-Gaussian and $X$ is a sub-Gaussian random vector. The above
inequality also implies
\[
\E\exp \bigl\{t\bigl(u^{\prime}X\bigr)/\sqrt{u^{\prime}\Sigma u}
\bigr\} \leq\exp\bigl(t^2\sigma^2/2\bigr).
\]
By Lemma~5.5 in  \cite{Vershynin:11}, there
exists a constant $c_0$
(depends only on $\sigma^2$) such that $\sqrt{c_0}\llVert (u^{\prime
}X)/\allowbreak \sqrt{u^{\prime}\Sigma u}\rrVert _{\psi_2} \leq1$. By the
linearity of the sub-Gaussian norm, we have $c_0\|u^{\prime} X\|_{\psi
_2}^2 \leq u^{\prime}\Sigma u$ as desired.
\end{pf*}

\subsubsection{Proofs of Theorems \texorpdfstring{\protect\ref{t:frobenius}}{2.1} and \texorpdfstring{\protect\ref{t:operator}}{2.2}}\label{sec:sample:proof}
For our analysis, we write $\Sigma_n= \Sigma_n^{\ast} - \bar{X}\bar
{X}^{\prime}$, where $\Sigma_n^{\ast} = n^{-1}\sum_{i=1}^n
X_iX_i^{\prime}$. Then $\|\Sigma_n-\Sigma\|_F \leq\|\Sigma_n^{\ast
}-\Sigma\|_F + \|\bar{X}\|^2$ and $\|\Sigma_n-\Sigma\|_2 \leq\|
\Sigma_n^{\ast}-\Sigma\|_2 + \|\bar{X}\|^2$. Hence to derive the
upper bounds for $\|\Sigma_n-\Sigma\|_F^2$ and $\|\Sigma_n-\Sigma\|
_2^2$, we just need to obtain the upper bounds for $\|\Sigma_n^{\ast
}-\Sigma\|_F^2$, $\|\Sigma_n^{\ast}-\Sigma\|_2^2$ and $\|\bar{X}\|
^4$. Because of the fact that $\P(X+Y\geq c+d) \leq\P(X\geq c) + \P
(Y\geq d)$ for any two univariate random variables $X$ and $Y$ and arbitrary
numbers $c$ and $d$, to study the tail behaviors of $\|\Sigma_n-\Sigma\|
_F$ and $\|\Sigma_n-\Sigma\|_2$, we only need to study those of $\|
\Sigma_n^{\ast}-\Sigma\|_F$, $\|\Sigma_n^{\ast}-\Sigma\|_2$ and
$\|\bar{X}\|^2$. As a result, Theorem~\ref{t:frobenius} is proved by
combining Propositions \ref{p:average} and~\ref{p:frobenius}, and
Theorem~\ref{t:operator} is proved by combining Propositions \ref
{p:average} and \ref{p:operator}. Materials that are needed for
proving Propositions \ref{p:frobenius} and \ref{p:frobenius} are
provided in the next two subsections.

We begin with the study of $\bar{X}\bar{X}^{\prime}$. Since this is
a rank 1 matrix, we make use of the basic fact
$\|\bar{X}\bar{X}^{\prime} \|_F = \| \bar{X}\bar{X}^{\prime}\|_2
= \| \bar{X} \|^2$. The following proposition is instrumental in the
proofs of Propositions \ref{p:average} and \ref{p:frobenius}.

%
\begin{prop}
\label{p:sub}
Let Assumption~\ref{as1} hold. There exist two fixed positive constants
$C_{\ast},
c_{\ast}$ such that, if $|t|>c_{\ast}(4c_0^{-1}+1) \tr(\Sigma)$,
\[
\E\exp \biggl\{\frac{\|X\|^2-\tr(\Sigma)}{t} \biggr\} \leq\exp \biggl\{C_{\ast}
\biggl[\frac{(4c_0^{-1}+1)\tr(\Sigma)}{t} \biggr]^2 \biggr\}.
\]
\end{prop}

%
\begin{pf}
Let $\|\cdot\|_{\psi_1}$
be the sub-exponential norm of a sub-exponential random variable (see
Definition~5.13 of  \cite{Vershynin:11}). We have
%
\begin{eqnarray}\label{eq_exponential}
\bigl\llVert \|X\|^2-\tr(\Sigma)\bigr\rrVert _{\psi_1} &\leq&
\bigl\llVert \|X\| ^2\bigr\rrVert _{\psi_1} + \bigl\llVert \tr(
\Sigma)\bigr\rrVert _{\psi
_1}
\nonumber
\\
& \leq&2\Vert \|X\|\Vert _{\psi_2}^2 + \tr(\Sigma)
\\
&\leq&\tr(\Sigma) \bigl(4c_0^{-1} + 1\bigr).\nonumber
\end{eqnarray}
For the second inequality above, we used Lemma~5.14 of  \cite{Vershynin:11} and for the third inequality we used
Lemma~\ref{l:sub_norm}. Because $\|X\|^2 - \tr(\Sigma)$ is a
zero-mean sub-exponential random variable, by Lemma~5.15 of  \cite{Vershynin:11}, there exist two fixed constants
$C_{\ast}, c_{\ast}$
such that if $|t|\geq c_{\ast}\llVert \|X\|^2-\tr(\Sigma)\rrVert
_{\psi_1}$,
\[
\E\exp \biggl\{\frac{\|X\|^2-\tr(\Sigma)}{t} \biggr\} \leq\exp \biggl\{C_{\ast}
\frac{\llVert \|X\|^2-\tr(\Sigma)\rrVert _{\psi
_1}^2}{t^2} \biggr\}.
\]
Combining (\ref{eq_exponential}) with the above inequality, we obtain
the proposition.
\end{pf}

%
\begin{prop}
\label{p:average}
Let Assumption~\ref{as1} hold. For any $t\geq0$,
%
\begin{equation}
\label{ineq_exp1} \P \biggl\{\|\bar{X}\|^2 \geq\frac{1+c_1t}{n}\cdot
\tr(\Sigma ) \biggr\} \leq\exp(1-t),
\end{equation}
where $c_1 = \max \{\max (\sqrt{C_{\ast}},c_{\ast}
)(4c_0^{-1}+1), 2 \}$ is a constant. Furthermore,
\[
\E \bigl(\llVert \bar{X}\rrVert ^4 \bigr) \leq \bigl\{ 1+2
\bigl(c_1^2+c_1\bigr)\exp(1) \bigr\}
\frac{\tr(\Sigma)^2}{n^2}.
\]
\end{prop}

%
\begin{pf}
It is straightforward
to verify that $\sqrt{n}\bar{X}$ is sub-Gaussian and satisfies
Assumption~\ref{as1} with the same $c_0$.
Applying the Markov inequality to $\exp (n\|\bar{X}\|^2 )$
we obtain, for any $a>0,  x\geq c_{\ast}(4c_0^{-1}+1)\tr(\Sigma)$,
\begin{eqnarray*}
\P \bigl\{n\|\bar{X}\|^2-\tr(\Sigma) \geq a \bigr\} &\leq&\exp
\bigl(-at^{-1}\bigr) \E\exp \bigl\{x^{-1} \bigl[n\|\bar{X}
\|^2-\tr(\Sigma ) \bigr] \bigr\}
\\
&\leq&\exp\bigl(-ax^{-1}\bigr) \exp \biggl\{C_{\ast} \biggl[
\frac
{(4c_0^{-1}+1)\tr(\Sigma)}{x} \biggr]^2 \biggr\},
\end{eqnarray*}
where the last inequality holds by Proposition~\ref{p:sub}. By letting
$ x= c_1\tr(\Sigma)$ and $ a = tx $ we obtain (\ref{ineq_exp1}). The
expectation inequality is proved in the supplemental article (\cite
{Bunea:14b}).
\end{pf}

Next, we study $\Sigma_n^{\ast} - \Sigma$. Let $Z_i = X_i
X_i^{\prime} -\Sigma$. Then $\E(Z_i) = 0$ and $\Sigma_n^* - \Sigma
= n^{-1}\sum_{i=1}^n Z_i$. We begin by stating the bounds with respect
to the Frobenius norm.

%
\begin{prop}
\label{p:frobenius}
Let Assumption~\ref{as1} hold. For all $n\geq1$ and $t\geq0$:
%
\begin{equation}
\label{ineq_exp_2} \P \biggl\{\bigl\|\Sigma_n^{\ast}-\Sigma
\bigr\|_F \geq\frac{2c_1 [
\sqrt{2\exp(1)} + 8\sqrt{t} ]\cdot\tr(\Sigma)}{\sqrt {n}} \biggr\} \leq2\exp \bigl\{-\min (t, 2
\sqrt{nt} ) \bigr\},
\end{equation}
where $c_1$ is defined in Proposition~\textup{\ref{p:average}}. Furthermore,
\[
\E \bigl(\bigl\|\Sigma_n^{\ast}-\Sigma\bigr\|_F^2
\bigr)\leq \biggl[ \frac
{4c_1\tr(\Sigma)}{\sqrt{n}} \biggr]^2 c_2,
\]
where $c_2 =\exp(1)+ \int_0^{\infty} \exp \{-\frac{1}{64}\min
(t,16\sqrt{t} ) \}\,\d t$.
\end{prop}
%

%
\begin{pf}
By Theorem~\ref
{t:JN2}, the Frobenius norm is $2$-smooth on the space $\mathbb
{R}^{p\times p}$ of $p\times p$ real matrices. Hence by Proposition~\ref{p:frobenius_exp} and Theorem~\ref{t:JN1},
\[
\P \biggl\{\bigl\|\Sigma_n^{\ast}-\Sigma\bigr\|_F \geq
\frac{2c_1 [
\sqrt{2\exp(1)} + t ]\cdot\tr(\Sigma)}{\sqrt{n}} \biggr\} \leq2\exp \biggl\{-\frac{1}{64}\min
\bigl(t^2, 16t\sqrt{n}\bigr) \biggr\}.
\]
Inequality (\ref{ineq_exp_2}) follows by changing $t$ to $8\sqrt{t}$
in the above inequality. The expectation inequality is proved in the
supplemental article (\cite{Bunea:14b}).
\end{pf}

%

%
\begin{prop}
\label{p:operator}
Let Assumption~\ref{as1} hold. For all $n\geq1$ and $t\geq0$:
%
\begin{equation}
\label{ineq_exp_3} \P \biggl\{\bigl\|\Sigma_n^{\ast}-\Sigma
\bigr\|_2\geq c_3\cdot\|\Sigma\| _2\cdot\max
\biggl\{\sqrt{\frac{r_e(\Sigma)(t+\ln p)}{n}},\frac
{r_e(\Sigma)(t+\ln p)}{n} \biggr\} \biggr\}\leq
\exp(-t),
\end{equation}
where $c_3$ is a fixed constant that depends only on $c_0$. Furthermore,
\[
\E \bigl(\bigl\|\Sigma_n^{\ast}-\Sigma\bigr\|_2^2
\bigr)\leq5c_3^2 \cdot\| \Sigma\|_2^2
\cdot\max \biggl\{\frac{r_e(\Sigma)\cdot\ln p }{n}, \biggl(\frac{r_e(\Sigma)\cdot\ln p }{n}
\biggr)^2 \biggr\}.
\]
\end{prop}
%

%
\begin{pf}
Let $Z_i =
X_iX_i^{\prime}-\Sigma$, then $\E(Z_i)=0$. We derive that $\Sigma
_n^{\ast} = n^{-1}\sum_{i=1}^n X_iX_i^{\prime} = n^{-1}\times \sum_{i=1}^n
Z_i + \Sigma$ and hence $\|\Sigma_n^{\ast}-\Sigma\|_2 = \|
n^{-1}\sum_{i=1}^n Z_i\|_2$. With Proposition~\ref{p:operator_sub3},
the probability inequality (\ref{ineq_exp_3}) is proved by applying
Theorem~\ref{t:JN3}. The expectation inequality is proved in the
supplemental article (\cite{Bunea:14b}).
\end{pf}

\subsubsection{Supplemental materials for proving Proposition \texorpdfstring{\textup{\protect\ref{p:frobenius}}}{A.3}}

The proof of Proposition~\texorpdfstring{\ref{p:frobenius}}{A.3} consists in adapting
results in
 \cite{Juditsky:08} to our context and
verifying its hypotheses. For
completeness, we state these results below.

%
\begin{theorem}
\label{t:JN1}
Let $(E, |\!|\!| \cdot|\!|\!| )$ be $\kappa$-smooth with a norm $|\!|\!| \cdot
|\!|\!|$ on $E$. Let $\{Z_1, Z_2,\dots\}$ be $E$-valued, zero-mean and
independent. Assume that there exists a sequence of positive numbers $\{
\sigma_1,\sigma_2,\dots\}$ such that
$
\E \{\exp (\sigma_i^{-1}|\!|\!|Z_i|\!|\!| ) \} \leq
\exp(1), i\geq1$.
Then for all $n\geq1$ and $t\geq0$:
\begin{eqnarray*}
\P \Biggl\{\biggl|\!\biggl|\!\biggl|\frac{Z_1+\cdots+ Z_n}{n}\biggr|\!\biggr|\!\biggr| \geq\frac{\sqrt{\exp
(1)\kappa}+t}{n}\sqrt{\sum
_{i=1}^n \sigma_i^2}
\Biggr\}\leq2\exp \biggl\{-\frac{1}{64}\min \bigl(t^2, t
t_n^{\ast} \bigr) \biggr\},
\end{eqnarray*}
where $ t_n^{\ast} = 16\sqrt{\sum_{i=1}^n \sigma_i^2}/\max_{1\leq
i\leq n} \sigma_i$.
\end{theorem}

%
\begin{rem}
Theorem~\ref{t:JN1} is a special case of Theorem~4.1 in  \cite{Juditsky:08} and the definition of a $\kappa
$-smooth space is on page
3 therein.
\end{rem}
%

%
\begin{theorem}
\label{t:JN2}
Let $2\leq p< \infty$. The Schatten norm $\|Z\|_p =  \{\sum_j
[d_j(Z) ]^p \}^{1/p}$ on the space $\mathbb
{R}^{m\times n}$ of $m\times n$ real matrices, where $d_1(Z)\geq
d_2(Z)\geq\cdots$ are the singular values of $Z$, is $\kappa_p(m,
n)$-smooth with
\[
\kappa_p(m,n) = \min_{2\leq\rho<\infty, \rho\leq p} \bigl\{\max (2,\rho-1)
\bigr\} \bigl\{\min(m,n) \bigr\}^{2/\rho- 2/p}.
\]
\end{theorem}

%
\begin{rem}
Theorem~\ref{t:JN2} is Example~3.3 in  \cite
{Juditsky:08}. For $p =
2$ we have the Frobenius norm which is $\kappa$-smooth with $\kappa= 2$.
\end{rem}

%
\begin{prop}
\label{p:frobenius_exp}
Let $Z = XX^{\prime}-\Sigma$. Then $
\E \{\exp [t^{-1}\|Z_i\|_F ] \}\leq\exp(1)$, for any $t \geq2c_1\tr(\Sigma)$, where $c_1$ is defined in
Proposition~\textup{\ref{p:average}}.
\end{prop}

\begin{pf}
First we have $\|
Z \|_F = \|XX^{\prime}-\Sigma\|_F \leq\|XX^{\prime}\|_F + \|\Sigma
\|_F = \|X\|^2 + \|\Sigma\|_F$. It is easy to show that $\|\Sigma\|_F
\leq\tr(\Sigma)$. Hence,
\begin{eqnarray*}
\E \bigl\{\exp \bigl[t^{-1}\|Z\|_F \bigr] \bigr\}&\leq&
\exp \bigl\{t^{-1} \bigl[\|\Sigma\|_F+\tr(\Sigma) \bigr]
\bigr\}\E \bigl\{ \exp \bigl[t^{-1}\bigl(\|X\|^2-\tr(\Sigma)
\bigr) \bigr] \bigr\}
\\
&\leq&\exp \bigl\{2t^{-1}\tr(\Sigma) \bigr\}\exp \biggl\{C \biggl[
\frac{(4c_0^{-1}+1)\tr(\Sigma)}{t} \biggr]^2 \biggr\}
\\
&\leq&\exp(1)
\end{eqnarray*}
as desired if $t>2c_1\tr(\Sigma)$. In the above derivation, we used
Proposition~\ref{p:sub}.
\end{pf}

\subsubsection{Supplemental materials for proving Proposition \texorpdfstring{\textup{\protect\ref{p:operator}}}{A.4}}

To derive the set of bounds on $\|\Sigma_n - \Sigma\|_2$ presented in
Proposition~\ref{p:operator}, we will appeal to the following result,
which is adapted from Theorem~6.2 in  \cite{Tropp:11}.

%
\begin{theorem}
\label{t:JN3}
Let $\{Z_i, i=1,\dots, n\}$ be a sequence of independent and
identically distributed symmetric matrices of dimension $p$. Assume
that there exist positive quantities $R$ and $\sigma$ such that
%
\begin{equation}
\label{cond} \E(Z_i) = 0 \quad \mbox{and}\quad  \bigl\llVert \mathrm{E}
\bigl(Z_i^d\bigr)\bigr\rrVert _2 \leq
\frac{d!}{2}\cdot R^{d-2}\sigma^2 \qquad \mbox{for } d =2,3,
\ldots.
\end{equation}
Then for all $t\geq0$, with probability at least $1 - \exp(-t)$,
\[
\biggl\llVert \frac{Z_1+\cdots+Z_n}{n}\biggr\rrVert _2< 3\cdot\max \biggl
\{ \sigma\sqrt{\frac{t+\ln p}{n}}, R\frac{t+\ln p}{n} \biggr\}.
\]
\end{theorem}

The proof of Proposition~\ref{p:operator} consists in the non-trivial
verification of condition (\ref{cond}). We do this in the following
proposition and two lemmas.

%
\begin{prop}
\label{p:operator_sub3}
Let Assumption~\ref{as1} hold, and define $Z = XX^{\prime} - \Sigma$, where
$\Sigma$ is the covariance matrix of $X$. Let $\tilde{c}_1 = \sup_{d\geq1} \exp(-d)d^d/d!$, $\tilde{c}_2 = \tilde{c}_1c_0^2\exp(-1)
+ \tilde{c}_1\exp(-1)/4 + 3$ and $\tilde{c}_3 = \max \{4\exp
(1)/c_0,1 \}$. If we let $R = 2\tilde{c}_3\cdot\tr(\Sigma)$
and $\sigma^2 = \tilde{c}_2 \tilde{c}_3^2\cdot\tr(\Sigma)\cdot\|
\Sigma\|_2$, then
\[
\bigl\llVert \E\bigl(Z^d\bigr)\bigr\rrVert _2 \leq
\frac{d!}{2}\cdot R^{d-2}\sigma^2 \qquad \mbox{for } d =2,3,
\ldots.
\]
\end{prop}

%

\begin{lem}
\label{l:sub4}
Suppose $A, B\in\mathbb{R}^{p\times p}$ are two positive
semi-definite matrices. Let $ODO^{\prime}$ be an eigendecomposition of
$A - B$ with $D = \diag(\lambda_1, \dots, \lambda_p)$. Let $D^{+}
= \diag(|\lambda_1|, \dots, |\lambda_p|)$. Then $OD^{+}O^{\prime}
\leq A + 2\|B\|_2 \cdot I_p$, where the notation \textup{``}$\leq$\textup{''} was used to
compare two matrices and for two matrices $E_1$ and $E_2$, $E_1\leq
E_2$ implies $E_2-E_1$ is psd.
\end{lem}


\begin{lem}
\label{l:sub5}
Suppose $A, B\in\mathbb{R}^{p\times p}$ are two positive
semi-definite matrices. Fix $u\in\mathbb{R}^p$. For an arbitrary
positive integer $d$,
\[
u^{\prime} (A-B)^d u \leq\|A-B\|_2^{d-1}
\bigl\{ u^{\prime}\bigl(A+ 2\|B\| _2\cdot I_p\bigr)u \bigr\}.
\]
\end{lem}


\subsubsection{Proof of Theorem \texorpdfstring{\protect\ref{t:trace}}{2.3}}\vspace*{-12pt}

\begin{pf*}{Proof of Theorem~\ref{t:trace}} Observe that $\tr(\Sigma
_n) = \tr(\Sigma_n^{\ast}) + \|\bar{X}\|^2$. With Proposition~\ref
{p:average}, it suffices to show that
\[
\P \bigl\{\bigl\llvert \tr\bigl(\Sigma_n^{\ast}\bigr)-\tr(
\Sigma)\bigr\rrvert \geq 2c_1\sqrt{t/n} \cdot\tr(\Sigma) \bigr\}\leq2
\exp(-t)
\]
for any $t\geq0$.
By the Markov inequality, if $nx \geq c_1\tr(\Sigma)$,
\begin{eqnarray*}
\P \bigl\{\tr\bigl(\Sigma_n^{\ast}\bigr)-\tr(\Sigma)\geq a
\bigr\}&\leq& \exp\bigl(-ax^{-1}\bigr) \E\exp \bigl\{x^{-1}
\bigl[\tr\bigl(\Sigma_n^{\ast
}\bigr)-\tr(\Sigma) \bigr] \bigr
\}
\\
&\leq&\exp\bigl(-ax^{-1}\bigr) \bigl\{\E\exp \bigl\{n^{-1}x^{-1}
\bigl[\|X\| ^2-\tr(\Sigma) \bigr] \bigr\} \bigr\}^n
\\
&\leq&\exp\bigl(-ax^{-1}\bigr)\exp \biggl\{C^{\ast} \biggl[
\frac
{(4c_0^{-1}+1)\tr(\Sigma)}{\sqrt{n}x} \biggr]^2 \biggr\},
\end{eqnarray*}
where in the last inequality we used Proposition~\ref{p:sub}. By
letting $x = c_1\tr(\Sigma)/\sqrt{nt}$ and $a = 2c_1\tr(\Sigma
)\cdot\sqrt{t/n}$ we obtain from the above inequality that
\[
\P \bigl\{\tr\bigl(\Sigma_n^{\ast}\bigr)-\tr(\Sigma)
\geq2c_1\sqrt{t/n} \cdot\tr(\Sigma) \bigr\} \leq\exp(-t).
\]
With a similar argument, we can obtain
\[
\P \bigl\{\tr\bigl(\Sigma_n^{\ast}\bigr)-\tr(\Sigma)\leq-
2c_1\sqrt{t/n} \cdot\tr(\Sigma) \bigr\} \leq\exp(-t)
\]
which completes the proof.
\end{pf*}

\subsubsection{Bounds on \texorpdfstring{$r_e(\Sigma_n)$}{$r_e(Sigma_n)$}}\label{sec:thm_re}

%
\begin{theorem}
\label{t:effective}
Suppose $X$ is a random vector that satisfies Assumption~\ref{as1}. Let $n>1$.
If $\Sigma\in\mathcal{P}_1(\epsilon)$, then with probability $1-11n^{-1}$,
\[
\biggl\llvert \frac{r_e(\Sigma_n)}{r_e(\Sigma)}-1\biggr\rrvert \lesssim\max \biggl\{\sqrt{
\frac{r_e(\Sigma)\cdot\ln pn }{2n}}, \frac
{r_e(\Sigma)\cdot\ln pn }{n} \biggr\}.
\]
If $\Sigma\in\mathcal{P}_2(\epsilon)$, then with probability $1-11n^{-1}$,
\[
\biggl\llvert \frac{r_e(\Sigma_n)}{r_e(\Sigma)}-1\biggr\rrvert \lesssim\frac
{r_e(\Sigma)\cdot\ln n }{n}.
\]
\end{theorem}

%

\subsection{Technical proofs of Section \texorpdfstring{\protect\ref{sec:detectable}}{3}}\label
{sec:reduced:proof}\vspace*{-12pt}

\begin{pf*}{Proof of Theorem~\ref{t:correct}} The proof follows from
arguments similar to those used in Theorem~2 of
\cite{Bunea:11b}. We
sketch it here for completeness.
Note that $\widehat{s}(\widetilde{\tau}) = s$ is equivalent to
$\widehat{\lambda}_s \geq\widetilde{\tau}$ and $\widehat{\lambda
}_{s+1} < \widetilde{\tau}$, or equivalently, $\lambda_s -\widehat
{\lambda}_s\leq\lambda_s - \widetilde{\tau}$ and $\widehat
{\lambda}_{s+1}-\lambda_{s+1} \leq\widetilde{\tau}-\lambda
_{s+1}$. By Weyl's theorem,\vspace*{1pt} Theorem~\ref{t:frobenius} and Theorem~\ref
{t:operator}, with probability larger than $1-5n^{-1}$, $\llvert \widehat{\lambda}_k - \lambda_k\rrvert  \leq\llVert \Sigma_n -
\Sigma\rrVert _2 \leq\eta_j$, for all $k$. Therefore, with (\ref
{tau}), it suffices to have $\lambda_s -\tau_1 \geq\eta_j$ and
$\tau_2-\lambda_{s+1}\geq\eta_j$, which is (\ref{index}).
\end{pf*}

\begin{pf*}{Proof of Theorem~\ref{t:detect}}
The proof is an application of Theorem~\ref{t:correct} with $\tau_1 =
2(1+\epsilon_1)\eta_j/C_j$ and $\tau_2 = 2C_j(1-\epsilon_1)\eta
_j$, and we just need to verify inequality (\ref{tau}) for
appropriately chosen $\delta$. By Theorem~\ref{t:trace}, with
probability $1-5n^{-1}$,
$
\llvert \tr(\Sigma_n)-\tr(\Sigma)\rrvert  \leq\epsilon_1 \tr
(\Sigma)$.
Let $\epsilon_2 = (1+c_1+c_3)\sqrt{\epsilon}$. For $\Sigma\in
\mathcal{P}_1(\epsilon)$, by Theorem~\ref{t:operator}, with
probability at least $1-4n^{-1}$,
$
\llVert \Sigma_n-\Sigma\rrVert _2 \leq\epsilon_2 \|\Sigma\|_2$.
Therefore, it is easy to show that, for $\Sigma\in\mathcal
{P}_1(\epsilon)$, with probability at least $1-6n^{-1}$,
$
\sqrt{(1-\epsilon_1)(1-\epsilon_2)}\eta_1 \leq\widetilde{\eta}_1
\leq\sqrt{(1+\epsilon_1)(1+\epsilon_2)}\eta_1$,
and
$
0.9 (1-\epsilon_1) \eta_1 \leq\widetilde{\eta}_1 \leq(1+\epsilon
_1)\eta_1/0.9
$
with the assumption that $\epsilon_2 \leq0.19$. For $\Sigma\in
\mathcal{P}_2(\epsilon)$, with probability at least $1-5n^{-1}$,
$
(1-\epsilon_1)\eta_2 \leq\widetilde{\eta}_2 \leq(1+\epsilon
_1)\eta_2
$.
\end{pf*}

\begin{pf*}{Proof of Theorem~\ref{t:detect_fda}}
The theorem is proved by combining Theorem~\ref{t:correct} and the
probability inequality $\P \{(1-\epsilon_1)\eta_2 \leq
\widetilde{\eta}_2 \leq(1+\epsilon_1)\eta_2 \}\geq1-
5n^{-1}$.
\end{pf*}

\subsection{Technical proofs of Section \texorpdfstring{\protect\ref{sec:eigen}}{4}}\label
{sec:eigen:proof}\vspace*{-12pt}
\begin{pf*}{Proof of Theorem~\ref{t:eigenvector_detect}} Note first that
\[
\min_{\lambda\in \mathit{EG}(\Sigma), \lambda\neq\lambda_k}|\lambda -\lambda_k| = \min(
\lambda_{k-1}-\lambda_k, \lambda_k - \lambda
_{k+1}),
\]
where we let $\lambda_0 = +\infty$ and $\lambda_{p+1} = 0$. By
Weyl's theorem and the results in Section~\ref{sec:eigen}, it is easy
to show that
\[
\min(\lambda_{k-1}-\lambda_k, \lambda_k -
\lambda_{k+1})\geq\min (\widehat{\lambda}_{k-1}-\widehat{
\lambda}_k, \widehat {\lambda}_k - \widehat{
\lambda}_{k+1} )- 2\eta_{\min},
\]
with probability larger than $1-5n^{-1}$. Because with probability
larger than $1-6n^{-1}$, $\widetilde{\eta}_j\geq C_j(1-\epsilon_1)
\eta_{\min}$,
the assumption (\ref{sample_spectrum}) in the theorem implies with
probability larger than $1-11n^{-1}$,
\[
\min(\lambda_{k-1}-\lambda_k, \lambda_k -
\lambda_{k+1})\geq3\eta _{\min}/\alpha,
\]
and the theorem holds by Proposition~\ref{p:eigenvector}.
\end{pf*}

\begin{pf*}{Proof of Theorem~\ref{t:eigen_detect}} Note that with
probability larger than $1-6n^{-1}$, $\widetilde{\eta}_{\mathrm{ev}} \geq
C_{1\lambda} (\frac{3\eta_2}{C_{3\lambda}\alpha}
)^{\beta_1/\beta_3} + \eta_2$. It follows that with probability
larger than $1-11n^{-1}$, $\lambda_k \geq\widehat{\lambda}_k -\eta
_2 \geq C_{1\lambda} (\frac{3\eta_2}{C_{3\lambda}\alpha
} )^{\beta_1/\beta_3}$, for all $k\leq\widetilde{K}_{\mathrm{ev}}$. By
Assumption~\ref{as2}, we derive that $k \leq (\frac{3\eta
_2}{C_{3\lambda}\alpha} )^{-1/\beta_3}$ and $\lambda_k -
\lambda_{k+1} \geq C_{3\lambda} k^{-\beta_3} \geq3\eta_2/\alpha$.
Therefore by Proposition~\ref{p:eigenvector}, with probability larger
than $1-11n^{-1}$, for all $k\leq\widetilde{K}_{\mathrm{ev}}$,
\[
\llVert \widehat{\boldsymbol{\psi}}_k- \boldsymbol{
\psi}_k\rrVert \leq\frac{\eta_2}{3\eta_2/\alpha} + \frac{6\eta_2^2}{9\eta
_2^2/\alpha^2}\leq\alpha,
\]
and
\[
\biggl\llvert \frac{\widehat{\lambda}_k}{\lambda_k} - 1\biggr\rrvert \leq\frac
{\eta_2}{\lambda_k} \leq
\frac{\eta_2}{\lambda_k - \lambda_{k+1}} \leq\frac{\alpha}{3}.
\]
\upqed
\end{pf*}

\subsection{Technical proofs of Section \texorpdfstring{\protect\ref{sec:fpca}}{5}}

\subsubsection{Proof of Proposition \texorpdfstring{\protect\ref{prop3}}{5.1}}\label
{sec:finite:proof}\vspace*{-12pt}
\begin{pf*}{Proof of Proposition~\ref{prop3}} First, notice that $\rho
_0$ is the integral $\int K(t,t)\,\mathrm{d}t$, while $\tr(\bK)=m^{-1}\sum_{j=1}^m K(t_j,t_j)$ is a finite approximation to the integral. Hence,
equality (\ref{trace})\vspace*{1pt} can be easily proved because of Assumption \ref{asD}.

To prove (\ref{eigenvalue}) and (\ref{eigenvector}), we need some
initial derivations.
By Assumptions \ref{asD}, \ref{asE} and \ref{asF}, we have
%
\begin{equation}
\label{eq0} \bigl\llvert \boldsymbol{\phi}_{k_1}^{\prime}
\boldsymbol{\phi}_{k_2} - \delta_{k_1,k_2}\bigr\rrvert \leq
C_{7\lambda}\max(k_1,k_2)^{\gamma_1}/m
\end{equation}
for all $k_1$ and $k_2$. Here $C_{7\lambda}$ is a fixed constant that
depends only on $C_{6\lambda}$ in Assumption \ref{asD} and $\delta_{k_1,k_2}$
equals 1 if $k_1=k_2$ and $0$ otherwise. Let $\lceil x \rceil$ be the
smallest integer that is no smaller than $x$. Define $N = \lceil
m^{1/(\beta_1+\gamma_1)}\rceil< m$. Let $A = [\boldsymbol{\phi
}_1,\dots, \boldsymbol{\phi}_N]$ be an $m\times N$ matrix and let $D
= \diag(\lambda_1,\dots, \lambda_N)$. It follows that
\[
\bK= \sum_{k} \lambda_k \boldsymbol{
\phi}_k \boldsymbol{\phi }_k^{\prime} = A
DA^{\prime} + \sum_{k>N} \lambda_k
\boldsymbol {\phi}_k \boldsymbol{\phi}_k^{\prime},
\]
and hence
\[
\bigl\llVert \bK- A DA^{\prime}\bigr\rrVert _F = \biggl
\llVert \sum_{k>N} \lambda _k
\boldsymbol{\phi}_k \boldsymbol{\phi}_k^{\prime}
\biggr\rrVert _F \leq\sum_{k>N}
\lambda_k\bigl\llVert \boldsymbol{\phi}_k\boldsymbol {
\phi}_k^{\prime}\bigr\rrVert _F = \sum
_{k>N} \lambda_k \|\boldsymbol {
\phi}_k\|^2.
\]
By Assumption \ref{asE}, $\lambda_k \leq C_{1\lambda}k^{-\beta_1}$. Hence,
\[
\sum_{k>N} \lambda_k \leq\int
_N^{\infty} C_{1\lambda}x^{-\beta
_1} \,\mathrm{d}x =
\frac{C_{1\lambda}}{1-\beta_1} x^{1-\beta_1}\biggl|_N^{\infty} =
\frac{C_{1\lambda}N^{1-\beta_1}}{\beta_1-1}.
\]
Combining the results above with (\ref{eq0}),
we obtain
%
\begin{equation}
\label{eq00} \bigl\llVert \bK- ADA^{\prime}\bigr\rrVert _F
= \sum_{k>N} \lambda_k\| \boldsymbol{
\phi}_k\|^2\leq C_{5\lambda}^2 \sum
_{k>N}\lambda_k \leq\frac{C_{5\lambda}^2C_{1\lambda}}{\beta_1-1}N^{1-\beta_1},
\end{equation}
where $C_{5\lambda}$ is an upper bound for all $\phi_k$ (see
Assumption \ref{asB}). Next, we study the term $ADA^{\prime}$. Consider a $QR$
decomposition of $A$, where $Q$ is an $ m\times N$ matrix with
orthonormal columns and $R$ is an $N\times N$ upper-triangular matrix.
Then $ADA^{\prime} = Q (RDR^{\prime} )Q^{\prime}$. Let
$Q$ and $R$ be given as in Lemma~\ref{lem5} below. We can further
derive for all $1\leq i, k\leq N$,
\[
\bigl\llvert R_{ik}^2 - \delta_{i,k}(1+r_i)^2
\bigr\rrvert \leq\frac
{5C_{7\lambda} k^{\gamma_1}}{m}\leq\frac{5C_{7\lambda} N^{\gamma_1}}{m}
\]
and for all $1\leq i, k, j\leq N$ with $i\neq j$,
\[
\llvert R_{ik}R_{jk}\rrvert \leq\frac{5C_{7\lambda} k^{\gamma
_1}}{m}\leq
\frac{5C_{7\lambda} N^{\gamma_1}}{m}.
\]
We let $\tilde{D} = RDR^{\prime}$ and compute $\tilde{d}_{ij}$
below. First,
\[
\tilde{d}_{ii} = \sum_k
\lambda_k R_{ik}^2 = \sum
_{1\leq k\leq N} \lambda_k \bigl\{R_{ik}^2
- \delta_{i,k}(1+r_i)^2 \bigr\} +\sum
_{1\leq k\leq N} \lambda_k \delta_{i,k}
(1+r_i)^2
\]
and hence
\[
\bigl\llvert \tilde{d}_{ii} - \lambda_i
(1+r_i)^2\bigr\rrvert \leq\sum
_{1\leq
k\leq N} \lambda_k \frac{5C_{7\lambda} N^{\gamma_1}}{m}=
\frac
{5C_{7\lambda}\rho_0 N^{\gamma_1}}{m}.
\]
Furthermore,
\begin{eqnarray*}
(\tilde{d}_{ii} -\lambda_i)^2 &\leq&\bigl(
\tilde{d}_{ii}- \lambda_i - 2\lambda_ir_i
- \lambda_ir_i^2\bigr)^2 +
\bigl(2\lambda_i r_i + \lambda _ir_i^2
\bigr)^2\\
&\leq&25\rho_0^2C_{7\lambda}^2N^{2\gamma_1}/m^2+
144\lambda _i^2 C_{7\lambda}^2N^{2+2\gamma_1}/m^2.
\end{eqnarray*}
Next for $i\neq j$,
\[
|\tilde{d}_{ij}| = \biggl\llvert \sum_{k}
\lambda_k R_{ik}R_{jk}\biggr\rrvert \leq
\frac{5\rho_0C_{7\lambda}N^{\gamma_1}}{m}.
\]
It follows that
\begin{eqnarray*}
\llVert \tilde{D} - D\rrVert _F^2&=&\sum
_{ij} (\tilde {d}_{ij}-
\lambda_i \delta_{ij} )^2
= \sum_i (\tilde{d}_{ii}-
\lambda_i )^2 +\sum_{i\neq
j}
\tilde{d}_{ij}^2
\\
&\leq& m^{-2}\sum_{i=1}^N \bigl
\{ 25\rho_0^2C_{7\lambda}^2 + 144
\lambda_i^2 C_{7\lambda}^2N^2
\bigr\}N^{2\gamma_1}+ m^{-2}\sum_{i\neq j}
25\rho_0^2 C_{7\lambda}^2N^{2\gamma_1}
\\
&\leq&169C_{7\lambda}^2\rho_0^2N^{2+2\gamma_1}/m^2,
\end{eqnarray*}
and hence
%
\begin{equation}
\label{orto} \bigl\llVert ADA^{\prime}- QDQ^{\prime}\bigr\rrVert
_F=\llVert \tilde{D} - D\rrVert _F\leq13C_{7\lambda}
\rho_0N^{1+\gamma_1}m^{-1}.
\end{equation}
Inequalities (\ref{eq00}) and (\ref{orto}) together lead to
%
\begin{equation}
\label{K_approx} \bigl\llVert \bK- QDQ^{\prime}\bigr\rrVert _F
\leq\frac{C_{5\lambda
}^2C_{1\lambda}N^{1-\beta_1}}{\beta_1-1} + \frac{13C_{7\lambda
}\rho_0N^{1+\gamma_1}}{m}.
\end{equation}

Now we are ready to prove (\ref{eigenvalue}) and (\ref{eigenvector}).
First, we invoke Weyl's theorem (\cite{Horn:85}, page
181), to
obtain, for each $k$,
%
\begin{eqnarray}\label
{eigenvalue_inequality}
\llvert \tilde{\lambda}_k - \lambda_k\rrvert 
&\leq& \bigl\llVert \bK-QDQ^{\prime}\bigr\rrVert _2
+ 1_{\{k>N\}}\lambda_k
\nonumber
\\
&\leq& \frac{C_{5\lambda}^2C_{1\lambda}N^{1-\beta_1}}{\beta_1-1} + \frac{13C_{7\lambda}\rho_0N^{1+\gamma_1}}{m} +C_{1\lambda}
N^{-\beta_1}
\\
&\leq& C_{8\lambda} m^{(1-\beta_1)/(\beta_1+\gamma_1)},\nonumber
\end{eqnarray}
where $C_{8\lambda}= C_{5\lambda}^2C_{1\lambda}/(\beta_1-1)
+C_{1\lambda} + 13C_{7\lambda}\rho_0$ is a fixed constant and recall
that $N = \lceil m^{1/(\beta_1+\gamma_1)}\rceil$.
Since the upper bound in the above derivation does not depend on $k$,
we obtain (\ref{eigenvalue}).

Finally, we prove (\ref{eigenvector}). As in Lemma~\ref{lem5} below,
we denote the columns of $Q$ by $\mathbf{v}_1,\dots, \mathbf{v}_N$.
Then for $1\leq k\leq N$, $\boldsymbol{\phi}_k = \sum_{j=1}^k R_{kj}
\mathbf{v}_j$. It follows that
%
\begin{eqnarray}
\label{eigenvector_approx1} \|\boldsymbol{\phi}_k - \mathbf{v}_k \|
&\leq&\sum_{j=1}^k |R_{kj}-
\delta_{k,j}| \leq|r_k| + \sum_{j=1}^k
3C_{7\lambda} j^{\gamma_1}/m \nonumber\\[-8pt]\\[-8pt]
&\leq&7C_{7\lambda} k^{1+\gamma_1}/m
\leq7C_{7\lambda
} N^{1+\gamma_1}/m.\nonumber
\end{eqnarray}
Next by Lemma A.1 in  \cite{Kneip:01} (see also
inequality (A.6) of
 \cite{Kneip:11}), we obtain from (\ref
{K_approx}) that
%
\begin{equation}
\label{eigenvector_approx2} \|\boldsymbol{\psi}_k - \mathbf{v}_k\|
\leq\frac{C_{8\lambda}
m^{(1-\beta)/(\beta+\gamma_1)}}{\min_{\lambda\in \mathit{EG}(\mathcal{K}),
\lambda\neq\lambda_k} |\lambda-\lambda_k| } +6 \biggl\{\frac{C_{8\lambda} m^{(1-\beta)/(\beta+\gamma_1)}}{\min_{\lambda\in \mathit{EG}(\mathcal{K}), \lambda\neq\lambda_k} |\lambda
-\lambda_k| } \biggr\}^2.
\end{equation}
Inequalities (\ref{eigenvector_approx1}) and (\ref
{eigenvector_approx2}) together gives (\ref{eigenvector}) which
completes the proof.
\end{pf*}

%
%
\begin{lem}
\label{lem5}
Suppose the assumptions in Proposition~\ref{prop3} hold. Let $A =
[\boldsymbol{\phi}_1,\dots, \boldsymbol{\phi}_N]$ be an $m\times
N$ matrix. Let $(Q, R)$ be a QR decomposition of $A$ where $Q$ is an $
m\times N$ matrix with orthonormal columns and $R$ is an $N\times N$
upper-triangular matrix. Denote the $(k,j)$th element of $R$ by
$R_{kj}$. Let $N$ be a positive integer such that $12C_{7\lambda}
N^{1+\gamma_1} \leq m$ where $C_{7\lambda}$ is the constant as in
inequality (\textup{\ref{eq0}}).
If $A$ has full rank, then there exists a pair of $Q$ and $R$ such that
if $k>j$, $R_{kj} = 0$ and if $k\leq j$,
\begin{eqnarray*}
\llvert R_{kj}-\delta_{k,j} -\delta_{k,j}
r_k \rrvert \leq 3C_{7\lambda}j^{\gamma_1}/m,
\end{eqnarray*}
where $r_k$ is defined in such a way that for all $k\leq N$
\[
|r_k| \leq4C_{7\lambda}k^{1+\gamma_1}/m.
\]
\end{lem}

\section*{Acknowledgements}
We thank David Sinclair for correcting an error in the paper.
Florentina Bunea was partially supported by Grant Number DMS 10007444
from National Science Foundation. Luo Xiao was partially supported by
Grant Number DMS 10007444 from National Science Foundation, Grant
Number R01EB012547 from the National Institute of Biomedical Imaging
And Bioengineering and Grant Number R01NS060910 from the National
Institute of Neurological Disorders and Stroke.

\end{appendix}


\begin{supplement}
\stitle{Supplement to ``On the sample covariance matrix estimator of reduced effective
rank population matrices, with applications to fPCA''\\}
\slink[doi]{10.3150/14-BEJ602SUPP} 
\sdatatype{.pdf}
\sfilename{BEJ602\_supp.pdf}
\sdescription{We provide proofs of all the lemmas, propositions and
theorems stated, but not proved, in the \hyperref[sec:proof]{Appendix} of
the main paper.}
\end{supplement}


\printhistory


\begin{thebibliography}{27}


\bibitem{Anderson:63}
\begin{barticle}[mr]
\bauthor{\bsnm{Anderson},~\bfnm{T.~W.}\binits{T.W.}}
(\byear{1963}).
\btitle{Asymptotic theory for principal component analysis}.
\bjournal{Ann. Math. Statist.}
\bvolume{34}
\bpages{122--148}.
\bid{issn={0003-4851}, mr={0145620}}
\end{barticle}
\bptok{imsref}%
\endbibitem

\bibitem{baing}
\begin{barticle}[mr]
\bauthor{\bsnm{Bai},~\bfnm{Jushan}\binits{J.}} \AND
\bauthor{\bsnm{Ng},~\bfnm{Serena}\binits{S.}}
(\byear{2002}).
\btitle{Determining the number of factors in approximate factor models}.
\bjournal{Econometrica}
\bvolume{70}
\bpages{191--221}.
\bid{doi={10.1111/1468-0262.00273}, issn={0012-9682}, mr={1926259}}
\end{barticle}
\bptok{imsref}%
\endbibitem

\bibitem{Baik:06}
\begin{barticle}[mr]
\bauthor{\bsnm{Baik},~\bfnm{Jinho}\binits{J.}} \AND
\bauthor{\bsnm{Silverstein},~\bfnm{Jack~W.}\binits{J.W.}}
(\byear{2006}).
\btitle{Eigenvalues of large sample covariance matrices of spiked population models}.
\bjournal{J. Multivariate Anal.}
\bvolume{97}
\bpages{1382--1408}.
\bid{doi={10.1016/j.jmva.2005.08.003}, issn={0047-259X}, mr={2279680}}
\end{barticle}
\bptok{imsref}%
\endbibitem

\bibitem{Benko:09}
\begin{barticle}[mr]
\bauthor{\bsnm{Benko},~\bfnm{Michal}\binits{M.}},
\bauthor{\bsnm{H{\"a}rdle},~\bfnm{Wolfgang}\binits{W.}} \AND
\bauthor{\bsnm{Kneip},~\bfnm{Alois}\binits{A.}}
(\byear{2009}).
\btitle{Common functional principal components}.
\bjournal{Ann. Statist.}
\bvolume{37}
\bpages{1--34}.
\bid{doi={10.1214/07-AOS516}, issn={0090-5364}, mr={2488343}}
\end{barticle}
\bptok{imsref}%
\endbibitem

\bibitem{Bickel:08b}
\begin{barticle}[mr]
\bauthor{\bsnm{Bickel},~\bfnm{Peter~J.}\binits{P.J.}} \AND
\bauthor{\bsnm{Levina},~\bfnm{Elizaveta}\binits{E.}}
(\byear{2008}).
\btitle{Covariance regularization by thresholding}.
\bjournal{Ann. Statist.}
\bvolume{36}
\bpages{2577--2604}.
\bid{doi={10.1214/08-AOS600}, issn={0090-5364}, mr={2485008}}
\end{barticle}
\bptok{imsref}%
\endbibitem

\bibitem{Bickel:08a}
\begin{barticle}[mr]
\bauthor{\bsnm{Bickel},~\bfnm{Peter~J.}\binits{P.J.}} \AND
\bauthor{\bsnm{Levina},~\bfnm{Elizaveta}\binits{E.}}
(\byear{2008}).
\btitle{Regularized estimation of large covariance matrices}.
\bjournal{Ann. Statist.}
\bvolume{36}
\bpages{199--227}.
\bid{doi={10.1214/009053607000000758}, issn={0090-5364}, mr={2387969}}
\end{barticle}
\bptok{imsref}%
\endbibitem

\bibitem{Bunea:11b}
\begin{barticle}[mr]
\bauthor{\bsnm{Bunea},~\bfnm{Florentina}\binits{F.}},
\bauthor{\bsnm{She},~\bfnm{Yiyuan}\binits{Y.}} \AND
\bauthor{\bsnm{Wegkamp},~\bfnm{Marten~H.}\binits{M.H.}}
(\byear{2011}).
\btitle{Optimal selection of reduced rank estimators of high-dimensional matrices}.
\bjournal{Ann. Statist.}
\bvolume{39}
\bpages{1282--1309}.
\bid{doi={10.1214/11-AOS876}, issn={0090-5364}, mr={2816355}}
\end{barticle}
\bptok{imsref}%
\endbibitem

\bibitem{Bunea:14b}
\begin{bmisc}[auto:STB|2014/02/12|14:17:21]
\bauthor{\bsnm{Bunea},~\bfnm{F.}\binits{F.}} \AND
\bauthor{\bsnm{Xiao},~\bfnm{L.}\binits{L.}}
(\byear{2014}). \bhowpublished{Supplement to ``On the sample covariance matrix estimator of
reduced effective rank population matrices, with applications to fPCA.'' DOI:\doiurl{10.3150/14-BEJ602SUPP}.}
\bid{doi={10.1350/14-BEJ602SUPP}}
\end{bmisc}
\bptok{imsref}%
\endbibitem

\bibitem{Cai:11}
\begin{barticle}[mr]
\bauthor{\bsnm{Cai},~\bfnm{Tony}\binits{T.}} \AND
\bauthor{\bsnm{Liu},~\bfnm{Weidong}\binits{W.}}
(\byear{2011}).
\btitle{Adaptive thresholding for sparse covariance matrix estimation}.
\bjournal{J. Amer. Statist. Assoc.}
\bvolume{106}
\bpages{672--684}.
\bid{doi={10.1198/jasa.2011.tm10560}, issn={0162-1459}, mr={2847949}}
\end{barticle}
\bptok{imsref}%
\endbibitem

\bibitem{Cai:10}
\begin{barticle}[mr]
\bauthor{\bsnm{Cai},~\bfnm{T.~Tony}\binits{T.T.}},
\bauthor{\bsnm{Zhang},~\bfnm{Cun-Hui}\binits{C.-H.}} \AND
\bauthor{\bsnm{Zhou},~\bfnm{Harrison~H.}\binits{H.H.}}
(\byear{2010}).
\btitle{Optimal rates of convergence for covariance matrix estimation}.
\bjournal{Ann. Statist.}
\bvolume{38}
\bpages{2118--2144}.
\bid{doi={10.1214/09-AOS752}, issn={0090-5364}, mr={2676885}}
\end{barticle}
\bptok{imsref}%
\endbibitem

\bibitem{Chamberlain:83}
\begin{barticle}[mr]
\bauthor{\bsnm{Chamberlain},~\bfnm{Gary}\binits{G.}} \AND
\bauthor{\bsnm{Rothschild},~\bfnm{Michael}\binits{M.}}
(\byear{1983}).
\btitle{Arbitrage, factor structure, and mean-variance analysis on large asset markets}.
\bjournal{Econometrica}
\bvolume{51}
\bpages{1281--1304}.
\bid{doi={10.2307/1912275}, issn={0012-9682}, mr={0736050}}
\end{barticle}
\bptok{imsref}%
\endbibitem

\bibitem{Dauxois:82}
\begin{barticle}[mr]
\bauthor{\bsnm{Dauxois},~\bfnm{J.}\binits{J.}},
\bauthor{\bsnm{Pousse},~\bfnm{A.}\binits{A.}} \AND
\bauthor{\bsnm{Romain},~\bfnm{Y.}\binits{Y.}}
(\byear{1982}).
\btitle{Asymptotic theory for the principal component analysis of a vector random function: Some applications to statistical inference}.
\bjournal{J. Multivariate Anal.}
\bvolume{12}
\bpages{136--154}.
\bid{doi={10.1016/0047-259X(82)90088-4}, issn={0047-259X}, mr={0650934}}
\end{barticle}
\bptok{imsref}%
\endbibitem

\bibitem{Fan:13}
\begin{barticle}[mr]
\bauthor{\bsnm{Fan},~\bfnm{Jianqing}\binits{J.}},
\bauthor{\bsnm{Liao},~\bfnm{Yuan}\binits{Y.}} \AND
\bauthor{\bsnm{Mincheva},~\bfnm{Martina}\binits{M.}}
(\byear{2013}).
\btitle{Large covariance estimation by thresholding principal orthogonal complements}.
\bjournal{J. R. Stat. Soc. Ser. B Stat. Methodol.}
\bvolume{75}
\bpages{603--680}.
\bid{doi={10.1111/rssb.12016}, issn={1369-7412}, mr={3091653}}
\bptnote{check related}%
\end{barticle}
\bptok{imsref}%
\endbibitem

\bibitem{Hall:06}
\begin{barticle}[mr]
\bauthor{\bsnm{Hall},~\bfnm{Peter}\binits{P.}} \AND
\bauthor{\bsnm{Hosseini-Nasab},~\bfnm{Mohammad}\binits{M.}}
(\byear{2006}).
\btitle{On properties of functional principal components analysis}.
\bjournal{J.~R. Stat. Soc. Ser. B Stat. Methodol.}
\bvolume{68}
\bpages{109--126}.
\bid{doi={10.1111/j.1467-9868.2005.00535.x}, issn={1369-7412}, mr={2212577}}
\end{barticle}
\bptok{imsref}%
\endbibitem

\bibitem{Hall:06b}
\begin{barticle}[mr]
\bauthor{\bsnm{Hall},~\bfnm{Peter}\binits{P.}},
\bauthor{\bsnm{M{\"u}ller},~\bfnm{Hans-Georg}\binits{H.-G.}} \AND
\bauthor{\bsnm{Wang},~\bfnm{Jane-Ling}\binits{J.-L.}}
(\byear{2006}).
\btitle{Properties of principal component methods for functional and longitudinal data analysis}.
\bjournal{Ann. Statist.}
\bvolume{34}
\bpages{1493--1517}.
\bid{doi={10.1214/009053606000000272}, issn={0090-5364}, mr={2278365}}
\end{barticle}
\bptok{imsref}%
\endbibitem

\bibitem{Horn:85}
\begin{bbook}[mr]
\bauthor{\bsnm{Horn},~\bfnm{Roger~A.}\binits{R.A.}} \AND
\bauthor{\bsnm{Johnson},~\bfnm{Charles~R.}\binits{C.R.}}
(\byear{1985}).
\btitle{Matrix Analysis}.
\blocation{Cambridge}:
\bpublisher{Cambridge Univ. Press}.
\bid{mr={0832183}}
\end{bbook}
\bptok{imsref}%
\endbibitem

\bibitem{Johnstone:01}
\begin{barticle}[mr]
\bauthor{\bsnm{Johnstone},~\bfnm{Iain~M.}\binits{I.M.}}
(\byear{2001}).
\btitle{On the distribution of the largest eigenvalue in principal components analysis}.
\bjournal{Ann. Statist.}
\bvolume{29}
\bpages{295--327}.
\bid{doi={10.1214/aos/1009210544}, issn={0090-5364}, mr={1863961}}
\end{barticle}
\bptok{imsref}%
\endbibitem

\bibitem{Juditsky:08}
\begin{bmisc}[auto:STB|2014/02/12|14:17:21]
\bauthor{\bsnm{Juditsky},~\bfnm{A.~B.}\binits{A.B.}} \AND
\bauthor{\bsnm{Nemirovski},~\bfnm{A.~S.}\binits{A.S.}}
(\byear{2008}). \bhowpublished{Large deviations of vector-values martingales
in 2-smooth normed spaces. Preprint. Available at \arxivurl{arXiv:0809.0813}}.
\end{bmisc}
\bptok{imsref}%
\endbibitem

\bibitem{Kneip:11}
\begin{barticle}[mr]
\bauthor{\bsnm{Kneip},~\bfnm{Alois}\binits{A.}} \AND
\bauthor{\bsnm{Sarda},~\bfnm{Pascal}\binits{P.}}
(\byear{2011}).
\btitle{Factor models and variable selection in high-dimensional regression analysis}.
\bjournal{Ann. Statist.}
\bvolume{39}
\bpages{2410--2447}.
\bid{doi={10.1214/11-AOS905}, issn={0090-5364}, mr={2906873}}
\end{barticle}
\bptok{imsref}%
\endbibitem

\bibitem{Kneip:01}
\begin{barticle}[mr]
\bauthor{\bsnm{Kneip},~\bfnm{Alois}\binits{A.}} \AND
\bauthor{\bsnm{Utikal},~\bfnm{Klaus~J.}\binits{K.J.}}
(\byear{2001}).
\btitle{Inference for density families using functional principal component analysis}.
\bjournal{J. Amer. Statist. Assoc.}
\bvolume{96}
\bpages{519--542}.
\bid{doi={10.1198/016214501753168235}, issn={0162-1459}, mr={1946423}}
\bptnote{check related}%
\end{barticle}
\bptok{imsref}%
\endbibitem

\bibitem{Lounici:13}
\begin{bmisc}[auto:STB|2014/02/12|14:17:21]
\bauthor{\bsnm{Lounici},~\bfnm{K.}\binits{K.}}
(\byear{2013}). \bhowpublished{High-dimensional covariance matrix estimation with missing observations. \emph{Bernoulli}. To appear}.
\end{bmisc}
\bptok{imsref}%
\endbibitem

\bibitem{muirhead}
\begin{bbook}[mr]
\bauthor{\bsnm{Muirhead},~\bfnm{Robb~J.}\binits{R.J.}}
(\byear{1982}).
\btitle{Aspects of Multivariate Statistical Theory}.
\bseries{Wiley Series in Probability and Mathematical Statistics}.
\blocation{New York}:
\bpublisher{Wiley}.
\bid{mr={0652932}}
\bptnote{check year}%
\end{bbook}
\bptok{imsref}%
\endbibitem

\bibitem{Nadler}
\begin{barticle}[mr]
\bauthor{\bsnm{Nadler},~\bfnm{Boaz}\binits{B.}}
(\byear{2008}).
\btitle{Finite sample approximation results for principal component analysis: A matrix perturbation approach}.
\bjournal{Ann. Statist.}
\bvolume{36}
\bpages{2791--2817}.
\bid{doi={10.1214/08-AOS618}, issn={0090-5364}, mr={2485013}}
\end{barticle}
\bptok{imsref}%
\endbibitem

\bibitem{Tropp:11}
\begin{barticle}[mr]
\bauthor{\bsnm{Tropp},~\bfnm{Joel~A.}\binits{J.A.}}
(\byear{2012}).
\btitle{User-friendly tail bounds for sums of random matrices}.
\bjournal{Found. Comput. Math.}
\bvolume{12}
\bpages{389--434}.
\bid{doi={10.1007/s10208-011-9099-z}, issn={1615-3375}, mr={2946459}}
\bptnote{check year}%
\end{barticle}
\bptok{imsref}%
\endbibitem

\bibitem{Vershynin:11}
\begin{bincollection}[mr]
\bauthor{\bsnm{Vershynin},~\bfnm{Roman}\binits{R.}}
(\byear{2012}).
\btitle{Introduction to the non-asymptotic analysis of random matrices}.
In \bbooktitle{Compressed Sensing}
\bpages{210--268}.
\blocation{Cambridge}:
\bpublisher{Cambridge Univ. Press}.
\bid{mr={2963170}}
\bptnote{check year}%
\end{bincollection}
\bptok{imsref}%
\endbibitem

\bibitem{Vu:12}
\begin{barticle}[auto:STB|2014/02/12|14:17:21]
\bauthor{\bsnm{Vu},~\bfnm{V.}\binits{V.}} \AND
\bauthor{\bsnm{Lei},~\bfnm{J.}\binits{J.}}
(\byear{2012}. \btitle{Minimax rates of estimation for sparse PCA in high
dimensions}.
\bjournal{JMLR: Workshop and Conference Proceedings}
\bvolume{22}
\bpages{1278--1286}.
\end{barticle}
\bptok{imsref}%
\endbibitem

\bibitem{Yao:05}
\begin{barticle}[mr]
\bauthor{\bsnm{Yao},~\bfnm{Fang}\binits{F.}},
\bauthor{\bsnm{M{\"u}ller},~\bfnm{Hans-Georg}\binits{H.-G.}} \AND
\bauthor{\bsnm{Wang},~\bfnm{Jane-Ling}\binits{J.-L.}}
(\byear{2005}).
\btitle{Functional data analysis for sparse longitudinal data}.
\bjournal{J. Amer. Statist. Assoc.}
\bvolume{100}
\bpages{577--590}.
\bid{doi={10.1198/016214504000001745}, issn={0162-1459}, mr={2160561}}
\end{barticle}
\bptok{imsref}%
\endbibitem

\end{thebibliography}
\end{document}